\documentclass[12pt]{article}
\input{epsf.sty}
\usepackage{epsfig}
\usepackage{amsmath}
\usepackage{latexsym}
\usepackage{amssymb}
\usepackage{bm}
\newtheorem{theorem}{Theorem}
\newtheorem{lemma}{Lemma}

\newtheorem{proposition}{Proposition}

\newtheorem{definition}{Definition}

\newtheorem{remark}{Remark}

\makeatletter
\let\normalequation=\equation
\def\equation{\@ifnextchar[{\subequation}{\normalequation}}
\def\subequation[#1]#2{\@ifundefined{r@#1}%
  {\def\theequation{\bf ??#2}\@warning
    {Reference `#1' on page \thepage \space
     undefined}}{\edef\@tempa{\@nameuse{r@#1}}%
    \edef\theequation{\expandafter\@car\@tempa \@nil#2}}%
  \let\@currentlabel\theequation $$}
\makeatother
\begin{document}

\begin{titlepage}
\begin{center}

\begin{center}
{\large\bf Desynchronization of pulse-coupled oscillators with
delayed excitatory coupling}\footnote{This work is supported
 by  National Science Foundation of China 60374018, 60574044 and
 973 Program 2005CB523306. }
\\[0.2in]
\begin{center}
Wei Wu\footnote{He is with Lab. of Nonlinear Mathematics Science,
Institute of Mathematics, Fudan University, Shanghai, 200433, P. R.
China   Email: wuweifd@vip.sina.com}, Tianping Chen\footnote{He is
with Lab. of Nonlinear Mathematics Science, Institute of
Mathematics, Fudan University, Shanghai, 200433, P. R. China.\\
\indent ~~Corresponding author: Tianping Chen. Email:
tchen@fudan.edu.cn; }
\end{center}
\end{center}
\end{center}
\begin{abstract}
Collective behavior of pulse-coupled oscillators has been
investigated widely. As an example of pulse-coupled networks,
fireflies display many kinds of flashing patterns. Mirollo and
Strogatz (1990) proposed a pulse-coupled oscillator model to explain
the synchronization of South East Asian fireflies ({\itshape
Pteroptyx malaccae}). However, transmission delays were not
considered in their model. In fact, the presence of transmission
delays can lead to desychronization. In this paper, pulse-coupled
oscillator networks with delayed excitatory coupling are studied.
Our main result is that under reasonable assumptions, pulse-coupled
oscillator networks with delayed excitatory coupling can not achieve
complete synchronization, which can explain why another species of
fireflies ({\itshape Photinus pyralis}) rarely synchronizes
flashing. Finally, two numerical simulations are given. In the first
simulation, we illustrate that even if all the initial phases are
very close to each other, there could still be big variations in the
times to process the pulses in the pipeline. It implies that
asymptotical synchronization typically also cannot be achieved. In
the second simulation, we exhibit a phenomenon of clustering
synchronization.
\end{abstract}


\end{titlepage}

\pagestyle{plain} \pagenumbering{arabic}
\section{Introduction}\quad
Fireflies provide one of the most spectacular examples of
synchronization in nature
\cite{Buck1976,Buck1988,Hanson1978,Smith1935}. Through observation,
people discovered that there are several synchronization patterns in
different species of fireflies. Some species, such as {\itshape
Pteroptyx malaccae}, flash rhythmically in perfect unison
\cite{Buck1976}. Some other species, such as {\itshape Photinus
pyralis}, usually show clustering synchrony or wave sweeping
synchrony, instead of complete synchrony \cite{Buck1988}.

Many kinds of biological models have been studied for flashing
behavior of fireflies. A pioneering work for synchronization of
fireflies has been reported by Buck in \cite{Buck1988}. He
investigated a variety of firefly species including {\itshape
Pteroptyx malaccae}, {\itshape Pteroptyx cribellata} and {\itshape
Photinus pyralis}. He suggested two kinds of flash models,
phase-advance entrainment model and phase-delay entrainment model.
In the phase advance model displayed in Fig.1(a), an excitation
level of pace-maker in the firefly brain is enhanced by an
external light stimulation. Each time excitation reaches a
threshold level, a neural signal is transmitted into the light
organ in the abdomen and then a flashing light is produced. An
example of this model is {\itshape Photinus pyralis}, the whole
group of the species rarely synchronizes flashing \cite{Buck1988}.
Instead, ``wave'', ``chain'' or ``sweeping'' synchrony has been
reported in some species with this phase advance model. A cluster
of male fireflies flash and then the flash of light is dispersed
into neighbor males. In the phase delay model displayed in
Fig.1(b), an excitation potential of pace-maker is reset to the
basal level by a light stimulation and the potential increase is
restarted to reach the threshold level for flashing. {\itshape
Pteroptyx cribellata} is known as an example of the phase delay
model and it usually shows a complete synchronization.

\setlength{\unitlength}{1cm}
\begin{picture}(20,7.9)
\put(0,5){\line(1,0){8.1}} \put(0,7){\line(1,0){8.1}}
\put(0.3,5){\line(-1,4){0.1}} \put(0.3,5){\line(6,5){2.4}}
\put(3.2,5){\line(-1,4){0.5}} \put(3.2,5){\line(6,5){1.8}}
\put(5,6.5){\line(1,5){0.1}} \put(5.6,5){\line(-1,4){0.5}}
\put(5.6,5){\line(6,5){1.8}} \put(7.4,6.5){\line(1,5){0.1}}
\put(7.9,5.4){\line(-1,4){0.4}} \put(0.3,4.95){\line(0,-1){0.3}}
\put(3.2,4.95){\line(0,-1){0.3}} \put(5.6,4.95){\line(0,-1){0.3}}
\put(1.5,4.77){\vector(-1,0){1.2}} \put(1.5,4.77){\vector(1,0){1.7}}
\put(4.5,4.77){\vector(-1,0){1.3}} \put(4.5,4.77){\vector(1,0){1.1}}
\put(5,7.6){\vector(0,-1){0.5}} \put(7.4,7.6){\vector(0,-1){0.5}}
\put(-0.5,7.8){\makebox(0,0)[c]{\footnotesize (a)}}
\put(0,7.2){\makebox(0,0)[c]{\scriptsize \itshape threshold}}
\put(-0.3,4.85){\makebox(0,0)[c]{\scriptsize \itshape base}}
\put(1.75,4.5){\makebox(0,0)[c]{\scriptsize \itshape normal period}}
\put(4.45,4.5){\makebox(0,0)[c]{\scriptsize \itshape updated period}}
\put(5.1,7.8){\makebox(0,0)[c]{\scriptsize \itshape light}}
\put(7.5,7.8){\makebox(0,0)[c]{\scriptsize \itshape light}}

\put(0,1){\line(1,0){8.1}} \put(0,3){\line(1,0){8.1}}
\put(0.3,1){\line(-1,4){0.1}} \put(0.3,1){\line(6,5){2.4}}
\put(3.2,1){\line(-1,4){0.5}} \put(3.2,1){\line(6,5){0.6}}
\put(3.9,1){\line(-1,5){0.1}} \put(3.9,1){\line(6,5){2.4}}
\put(6.8,1){\line(-1,4){0.5}} \put(6.8,1){\line(6,5){0.6}}
\put(7.5,1){\line(-1,5){0.1}} \put(7.5,1){\line(6,5){0.6}}
\put(0.3,0.95){\line(0,-1){0.3}} \put(3.2,0.95){\line(0,-1){0.3}}
\put(6.8,0.95){\line(0,-1){0.3}}
\put(1.5,0.77){\vector(-1,0){1.2}}
\put(1.5,0.77){\vector(1,0){1.7}} \put(5,0.77){\vector(-1,0){1.8}}
\put(5,0.77){\vector(1,0){1.8}} \put(3.8,3.6){\vector(0,-1){0.5}}
\put(7.4,3.6){\vector(0,-1){0.5}}
\put(-0.5,3.8){\makebox(0,0)[c]{\footnotesize (b)}}
\put(0,3.2){\makebox(0,0)[c]{\scriptsize \itshape threshold}}
\put(-0.3,0.85){\makebox(0,0)[c]{\scriptsize \itshape base}}
\put(1.75,0.5){\makebox(0,0)[c]{\scriptsize \itshape normal period}}
\put(5.1,0.5){\makebox(0,0)[c]{\scriptsize \itshape updated period}}
\put(3.8,3.8){\makebox(0,0)[c]{\scriptsize \itshape light}}
\put(7.4,3.8){\makebox(0,0)[c]{\scriptsize \itshape light}}

\put(3.7,-0.1){\makebox(0,0)[c]{\footnotesize Fig.1: (a) Phase
advance model (b) Phase delay model.}}
\put(1,-0.6){\makebox(0,0)[c]{\footnotesize (Reprinted from
\cite{Buck1988})}}

\put(10.5,6.5){\line(1,0){5}} \put(10.5,1.5){\line(1,0){5}}
\put(10.5,1.5){\line(0,1){5}} \put(15.5,1.5){\line(0,1){5}}
\qbezier(10.5,1.5)(11,6)(15.5,6.5) {\thicklines
\put(13,5.75){\vector(4,3){0.01}}}
\put(10.15,1.25){\makebox(0,0)[c]{\scriptsize \itshape (0,0)}}
\put(15.5,1.25){\makebox(0,0)[c]{\scriptsize \itshape (1,0)}}
\put(10.15,6.5){\makebox(0,0)[c]{\scriptsize \itshape (0,1)}}
\put(13,1){\makebox(0,0)[c]{\footnotesize $\varphi$}}
\put(10,4){\makebox(0,0)[c]{\footnotesize $x$}}
\put(11.5,5.3){\makebox(0,0)[c]{\scriptsize $x=f(\varphi)$}}
\put(12.5,-0.1){\makebox(0,0)[c]{\footnotesize Fig.2: Graph of the
function $f$. The time-course}}
\put(11.7,-0.6){\makebox(0,0)[c]{\footnotesize of the oscillator
is given by $x=f(\varphi)$.}}
\end{picture}

Mirollo and Strogatz (1990) argued that Buck's model is not
appropriate for synchronization because he assumed a linearly
increasing potential toward threshold as shown in Fig.1. Inspired
by Peskin's model for self-synchronization of the cardiac
pacemaker \cite{Peskin1975}, Mirollo and Strogatz proposed a
pulse-coupled oscillator model with undelayed excitatory coupling
to explain the synchronization of huge congregations of South East
Asian fireflies ({\itshape Pteroptyx malaccae}) \cite{Mir1990}.
This model is a network of $N$ pulse-coupled oscillators. Each
oscillator is characterized by a state variable $x_i$ which is
assumed to increase toward a threshold at $x_i=1$ according to
$x_i=f(\varphi_i)$, where $f: [0,1]\rightarrow[0,1]$ is smooth,
monotonic increasing, and concave down, i.e., $f'>0$ and $f''<0$.
Here $\varphi_i\in[0,1]$ is a phase variable such that (i)
$\mathrm{d}\varphi_i/\mathrm{d}t=1/T$, where $T$ is the cycle
period, (ii) $\varphi_i=0$ when the $i$th oscillator is at its
lowest state $x_i=0$, and (iii) $\varphi_i=1$ at the end of the
cycle when the $i$th oscillator reaches the threshold $x_i=1$.
Therefore $f$ satisfies $f(0)=0$, $f(1)=1$. Fig.2 shows the graph
of a typical $f$. When $x_i$ reaches the threshold, the $i$th
oscillator ``fires'' and $x_i$ jumps back instantly to zero, after
which the cycle repeats. The oscillators are assumed to interact
by a simple form of pulse-coupling: when a given oscillator fires,
it pulls all the other oscillators up by an amount $\varepsilon$,
or pulls them up to firing, whichever is less. That is,
\begin{eqnarray}
x_i(t)=1\Rightarrow x_j(t^{+})=\min(1,x_j(t)+\varepsilon),\quad
\forall j\neq i.\label{UndelayedInteraction}
\end{eqnarray}
The main result in \cite{Mir1990} is that for all $N$ and for
almost all initial conditions, the system eventually becomes
completely synchronized. Here, the concept of complete
synchronization is defined as follows: if there exists a
$t_0\geq0$ such that
\begin{eqnarray}
x_i(t)=x_j(t),\quad\mbox{for all $t\geq t_0$ and all $i\neq j$}
\end{eqnarray}
then the pulse-coupled oscillator network is said to be completely
synchronized or, for simplicity, synchronized.

Since Mirollo and Strogatz's model was introduced, many results on
pulse-coupled networks with undelayed coupling, including
undelayed excitatory coupling and undelayed inhibitory coupling,
have been obtained
\cite{Kuramoto1991,Vanvreeswijk1993,Goel2002,Chen1994,Corral1995,Mathar1996}.
However, reaction and transmission delays are unavoidable in real
biological systems. For example, normally the transmission delays
of most fireflies from sensors to motor actions of flashing are
around 200ms (see \cite{Buck1988}). In fact, transmission delays
influence the performance of synchronization for both excitatory
and inhibitory coupling \cite{Nischwitz1995,Ernst1995,Ernst1998}.
In particular, it has been shown that for excitatory coupling, the
presence of transmission delays can lead to desychronization
\cite{Nischwitz1995,Ernst1995,Knoblauch2002,Coombes1997}. To the
best of our knowledge, this desynchronization was proved only for
the case of two pulse-coupled oscillators, while for the case of
$N>2$ pulse-coupled oscillators, it was revealed only in
simulations.

In this paper, pulse-coupled oscillator networks with delayed
excitatory coupling are studied. We retain two of Mirollo and
Strogatz's assumptions: the oscillators have identical dynamics, and
all oscillators are coupled to all the others. Our main contribution
is that we prove that $N\geq2$ pulse-coupled oscillators with
delayed excitatory coupling can not achieve complete
synchronization. This result can explain why {\itshape Photinus
pyralis} rarely synchronizes flashing, which is known as an example
of pulse-coupled oscillator networks with delayed excitatory
coupling \cite{Buck1988}. Finally, two numerical simulations are
given. In Simulation 1, we illustrate that even if all the initial
phases are very close to each other, there could still be big
variations in the times to process the pulses in the pipeline. It
implies that asymptotical synchronization typically also cannot be
achieved. In Simulation 2, we exhibit a phenomenon of clustering
synchronization.

The rest of the paper is organized as follows: In Section 2, we
describe the network model. In Section 3, some definitions and
notations are given. In Section 4, some lemmas are proved. In
Section 5, we prove the main theorem. Two numerical simulations are
given in Section 6. We conclude the paper in Section 7.

\section{Model}\quad
The network consists of $N\geq2$ pulse-coupled oscillators with
delayed excitatory coupling. As in \cite{Mir1990}, the oscillators
have identical dynamics, and all oscillators are coupled to all
the others. Each oscillator is characterized by a state variable
$x_i$ which is assumed to increase toward a threshold at $x_i=1$
according to $x_i=f(\varphi_i)$, where $f: [0,1]\rightarrow[0,1]$
is smooth, monotonic increasing, and concave down, i.e., $f'>0$
and $f''<0$. Here $\varphi_i\in[0,1]$ is a phase variable such
that (i) $\mathrm{d}\varphi_i/\mathrm{d}t=1/T$, where $T$ is the
cycle period, (ii) $\varphi_i=0$ when the $i$th oscillator is at
its lowest state $x_i=0$, and (iii) $\varphi_i=1$ at the end of
the cycle when the $i$th oscillator reaches the threshold $x_i=1$.
When $x_i$ reaches the threshold, the $i$th oscillator fires and
$x_i$ jumps back instantly to zero, after which the cycle repeats.
That is,
\begin{eqnarray}
x_i(t)=1\Rightarrow x_i(t^{+})=0.\label{self-variety}
\end{eqnarray}

The new features are that, (i) in order to make the discussion a
little easier, we let $T=1$, i.e.,
$\mathrm{d}\varphi_i/\mathrm{d}t=1$, and (ii) the transmission
delay $\tau$ is introduced into the model. Because of the presence
of the transmission delay, the oscillators interact by a new form:
when a given oscillator fires at time $t$, it emits a spike; after
a transmission delay $\tau$, the spike reaches all the other
oscillators at time $t+\tau$ and pulls them up by an amount
$\varepsilon$, or pulls them up to firing, whichever is less. That
is,
\begin{eqnarray}
x_i(t)=1\Rightarrow
x_j(t+\tau)=\min(1,x_j((t+\tau)^{-})+\varepsilon),\quad \forall
j\neq i,\label{DelayedInteraction}
\end{eqnarray}
instead of (\ref{UndelayedInteraction}). Of course, if there are $m$
oscillators firing simultaneously, then the state variables $x$ of
the $m$ oscillators are increased by an amount $(m-1)\varepsilon$
and the state variables $x$ of all the other oscillators are
increased by an amount $m\varepsilon$. By (\ref{self-variety}) and
(\ref{DelayedInteraction}), it is clear that if a spike reaches the
$i$th oscillator just at the moment that the $i$th oscillator
reaches the threshold $x_i=1$, then the state variable $x_i$ changes
as follows:
\begin{eqnarray*}
x_i(t)=\min(1,x_i(t^{-})+\varepsilon)=\min(1,1+\varepsilon)=1\Rightarrow
x_i(t^{+})=0.
\end{eqnarray*}

In addition, we make the following two assumptions:
\newcounter{mycounter}
\begin{list}
{{\upshape (A\arabic{mycounter})}\hfill}
{\setlength{\topsep}{0ex}
 \setlength{\parskip}{0ex}
 \setlength{\itemsep}{0.2ex}
 \setlength{\parsep}{0ex}
 \setlength{\leftmargin}{6ex}
 \setlength{\labelwidth}{4ex}
 \setlength{\labelsep}{1ex}
 \setlength{\itemindent}{-1ex}
 \usecounter{mycounter}}
\item The system is started at time $t=0$ with a set of initial
states $0<x_i(0)\leq1$, and there are no firings in time
$[-\tau,0)$.

\item The transmission delay $\tau$ and the coupling strength
$\varepsilon$ satisfy $f(2\tau)+N\varepsilon<1$.
\end{list}

The assumption (A1) about the initial conditions is quite common
for pulse-coupled networks with transmission delays
\cite{Nischwitz1995,Ernst1995,Knoblauch2002,Coombes1997}. The time
$t=0$ can be regarded as the moment at which fireflies assemble
and begin to flash. Since $x_i$ jumps back instantly to $0$ when
it reaches the threshold $1$, the state values $0$ and $1$ may be
considered to be the same state. Therefore, we assume that the
initial states satisfy $0<x_i(0)\leq1$ instead of $0\leq
x_i(0)\leq1$.

The assumption (A2) means that both $\tau$ and $\varepsilon$ are
relatively small. It is reasonable to fireflies with the phase
advance entrainment model, especially to {\itshape Photinus
pyralis}. The transmission delays of most fireflies from sensors
to motor actions of flashing are around 200ms, while the regular
endogenous period of {\itshape Photinus pyralis} is almost 6000ms
(see \cite{Buck1988}). Normally, excitatory activations generate
frequent spikes in a short period of time, because it can shorten
the period of oscillator by increasing the potential level (see
\cite{Kim2004}). It means that strong excitations result in
frequent flashing of light, which will exhaust out the energy of
the insect. Thus, a relatively small coupling strength for
excitatory action is used in many species of fireflies, which
follow the phase advance model. In the following sections, one can
see that under this assumption, the number of case distinctions
will be reduced.

With the monotonicity of $f$, the state variable $x_i$ and the
phase variable $\varphi_i$ are one-to-one correspondence.
Therefore, the synchronization of the state variables
$x_1,\ldots,x_N$ is equivalent to the synchronization of the phase
variables $\varphi_1,\ldots,\varphi_N$. In the following, instead
of investigating $x_i$, we investigate dynamical behaviors of
$\varphi_i$ directly.

From the above description, the phase variable processes the
following properties.
\begin{proposition}\quad 
The phase variable $\varphi_i$ satisfies:
\begin{list}
{{\upshape (\alph{mycounter})}\hfill}
{\setlength{\topsep}{1ex}
 \setlength{\parskip}{0ex}
 \setlength{\itemsep}{0.2ex}
 \setlength{\parsep}{0ex}
 \setlength{\leftmargin}{3.75ex}
 \setlength{\labelwidth}{3.5ex}
 \setlength{\labelsep}{0.25ex}
 \usecounter{mycounter}}
\item If $m$ spikes reach the $i$th oscillator at time $t$, then
$\varphi_i(t)=f^{-1}\big(\min[1,f(\varphi_i(t^{-}))+m\varepsilon]\big)$;\\
If no spikes reach the $i$th oscillator at time $t$, then
$\varphi_i(t)=\varphi_i(t^{-})$.

\item If $\varphi_i(t)=1$, then $\varphi_i(t^{+})=0$;\\
If $\varphi_i(t)<1$, then $\varphi_i(t^{+})=\varphi_i(t)$.

\item If no spikes reach the $i$th oscillator in time $(t_1,t_2)$
and the $i$th oscillator do not fire in time $(t_1,t_2)$, then
$\varphi_i(t_2^{-})=\varphi_i(t_1^{+})+(t_2-t_1)$.
\end{list}
\end{proposition}

\begin{remark}\quad 
The assumption {\upshape (A1)} and Proposition {\upshape 1(b)}
imply that $\varphi_i(t)\neq0$ for all $t\geq0$ and all $1\leq
i\leq N$. Thus, in our model, each phase variable $\varphi_i$
satisfies $\varphi_i(t)\in(0,1]$ for all $t\geq0$.
\end{remark}

\section{Definitions and notations}\quad
In this paper, the concept of complete synchronization is defined
as follows:

\begin{definition}\quad 
If there exists a $t_0\geq0$ such that $\varphi_i(t)=\varphi_j(t)$
for all $t\geq t_0$, then we say that oscillator $i$ and $j$ can
achieve complete synchronization or, for simplicity,
synchronization; and say that oscillator $i$ and $j$ have already
been completely synchronized at time $t_0$ or, for simplicity,
synchronized.
\end{definition}

\begin{definition}\quad 
If there exists a $t_0\geq0$ such that
$\varphi_1(t)=\varphi_2(t)=\cdots=\varphi_N(t)$ for all $t\geq
t_0$, then we say that the pulse-coupled oscillator network can
achieve complete synchronization or, for simplicity,
synchronization; and say that the pulse-coupled oscillator network
has already been completely synchronized at time $t_0$ or, for
simplicity, synchronized.
\end{definition}

Throughout the paper, the following notations will be used.
\begin{list}
{(\roman{mycounter})\hfill}
{\setlength{\topsep}{1ex}
 \setlength{\parskip}{0ex}
 \setlength{\itemsep}{0.2ex}
 \setlength{\parsep}{0ex}
 \setlength{\leftmargin}{4.5ex}
 \setlength{\labelwidth}{4ex}
 \setlength{\labelsep}{0.5ex}
 \usecounter{mycounter}}
\item $F_m(\theta)=f^{-1}(\min[1,f(\theta)+m\varepsilon])$,
where $m\in \mathbb{Z}^{+}$ and $0\leq\theta\leq1$. Here,
$\mathbb{Z}^{+}=\{z\in\mathbb{Z}|z\geq0\}$.

\item $T_i=\{t\in \mathbb{R}^{+}\big|$t is the time at which
the $i$th oscillator fires$\}$, $i=1,\ldots,N$. Here,
$\mathbb{R}^{+}=\{z\in\mathbb{R}|z\geq0\}$.

\item $D_i(t)=\{\eta=t'+\tau-t\bigm|t'\in
T_i\cap[t-\tau,t)\}$, $i=1,\ldots,N$, $t\geq0$.

\item $D_i(t^{+})=\{\eta=t'+\tau-t\bigm|t'\in
T_i\cap(t-\tau,t]\}$, $i=1,\ldots,N$, $t\geq0$, which is called
the right-limit of $D_i(t)$ at time $t$.
\end{list}

\begin{remark}\quad 
According to Proposition 1, if $m$ spikes reach the $i$th
oscillator at time $t$, then $\varphi_i(t)
=f^{-1}\big(\min[1,f(\varphi_i(t^{-}))+m\varepsilon]\big)=F_m(\varphi_i(t^{-}))$;
if no spikes reach the $i$th oscillator at time $t$, then
$\varphi_i(t)=\varphi_i(t^{-})=F_0(\varphi_i(t^{-}))$. Therefore,
for the convenience of later use, we introduced the notation
$F_m(\theta)$.
\end{remark}

\begin{remark}\quad 
The notation $T_i$ is called the firing-time set of the $i$th
oscillator. Using the firing-time set $T_i$, we defined two other
notations $D_i(t)$ and $D_i(t^{+})$. They can be used to calculate
the reaching time of the spikes which the $i$th oscillator has
emitted, but the other oscillators have not received. From the
definition of $D_i(t)$, one can see that for each $\eta\in D_i(t)$,
there must be a spike of the $i$th oscillator, which will reach the
other oscillators at time $t+\eta$. In particular, if $0\in D_i(t)$,
then there must be a spike of the $i$th oscillator, which reaches
the other oscillators at time $t$.
\end{remark}

\begin{remark}\quad 
The notations $D_i(t)$ and $D_i(t^{+})$ will play an important role
in our discussion. For the pulse-coupled networks without
transmission delays, instantaneous synchronization
{\upshape(}$\varphi_i(t_0)=\varphi_j(t_0)${\upshape)} shows that
oscillator $i$ and $j$ have been completely synchronized at time
$t_0$ {\upshape(}$\varphi_i(t)=\varphi_j(t)$ for all $t\geq
t_0${\upshape)}, because the coupling is all-to-all and the
oscillators have identical dynamics. But, for those with
transmission delays, this criterion fails {\upshape(}instantaneous
synchronization does not mean complete synchronization{\upshape)},
because of the presence of delays. Fig.3 shows an example in which
the criterion fails. One key role of $D_i(t)$ and $D_i(t^{+})$ in
this paper is to judge whether the pulse-coupled oscillator network
with delays has been completely synchronized at time $t_0$, when we
have $\varphi_1(t_0)=\ldots=\varphi_N(t_0)$.
\end{remark}

\begin{picture}(20,5.25)
\put(0.5,1){\line(1,0){4}}
\put(0.5,5){\line(1,0){4}}
\put(0.5,1){\line(0,1){4}}
\put(4.5,1){\line(0,1){4}}
\qbezier(0.5,1)(1.5,4)(4.5,5)
\put(0.15,0.75){\makebox(0,0)[c]{\scriptsize \itshape (0,0)}}
\put(4.5,0.75){\makebox(0,0)[c]{\scriptsize \itshape (1,0)}}
\put(0.15,5){\makebox(0,0)[c]{\scriptsize \itshape (0,1)}}
\put(3.75,4.7){\circle*{0.15}}
\put(3.75,4.4){\makebox(0,0)[c]{\scriptsize A}}
\put(4.42,4.95){\circle{0.15}}
\put(4.4,4.7){\makebox(0,0)[c]{\scriptsize B}}
\put(3.75,1){\line(0,1){0.2}}
\put(3.75,0.75){\makebox(0,0)[c]{\scriptsize $1-\phi$}}
\put(2.5,0){\makebox(0,0)[c]{\footnotesize (a)}}

\put(6.5,1){\line(1,0){4}}
\put(6.5,5){\line(1,0){4}}
\put(6.5,1){\line(0,1){4}}
\put(10.5,1){\line(0,1){4}}
\qbezier(6.5,1)(7.5,4)(10.5,5)
\put(9.75,4.7){\circle*{0.15}}
\put(9.75,4.4){\makebox(0,0)[c]{\scriptsize A}}
\put(6.55,1.05){\circle{0.15}}
\put(6.75,1.15){\makebox(0,0)[c]{\scriptsize B}}
\put(9.75,1){\line(0,1){0.2}}
\put(9.75,0.75){\makebox(0,0)[c]{\scriptsize $1-\phi$}}
\put(8.5,0){\makebox(0,0)[c]{\footnotesize (b)}}

\put(12.5,1){\line(1,0){4}}
\put(12.5,5){\line(1,0){4}}
\put(12.5,1){\line(0,1){4}}
\put(16.5,1){\line(0,1){4}}
\qbezier(12.5,1)(13.5,4)(16.5,5)
\put(16.42,4.95){\circle*{0.15}}
\put(16.4,4.7){\makebox(0,0)[c]{\scriptsize A}}
\put(13.25,2.6){\circle{0.15}}
\put(13.44,2.4){\makebox(0,0)[c]{\scriptsize B}}
\put(13.25,1){\line(0,1){0.2}}
\put(13.25,0.75){\makebox(0,0)[c]{\scriptsize $\phi$}}
\put(14.5,0){\makebox(0,0)[c]{\footnotesize (c)}}
\end{picture}

\begin{picture}(20,5.25)
\put(0.5,0.5){\line(1,0){4}}
\put(0.5,4.5){\line(1,0){4}}
\put(0.5,0.5){\line(0,1){4}}
\put(4.5,0.5){\line(0,1){4}}
\qbezier(0.5,0.5)(1.5,3.5)(4.5,4.5)
\put(0.55,0.55){\circle*{0.15}}
\put(0.8,0.65){\makebox(0,0)[c]{\scriptsize A}}
\put(1.25,2.05){\circle{0.15}}
\put(1.44,1.9){\makebox(0,0)[c]{\scriptsize B}}
\put(1.25,0.5){\line(0,1){0.2}}
\put(1.25,0.25){\makebox(0,0)[c]{\scriptsize $\phi$}}
\put(2.5,-0.5){\makebox(0,0)[c]{\footnotesize (d)}}

\put(6.5,0.5){\line(1,0){4}}
\put(6.5,4.5){\line(1,0){4}}
\put(6.5,0.5){\line(0,1){4}}
\put(10.5,0.5){\line(0,1){4}}
\qbezier(6.5,0.5)(7.5,3.5)(10.5,4.5)
\put(7.6,2.6){\circle*{0.15}}
\put(7.7,2.9){\makebox(0,0)[c]{\scriptsize A}}
\put(7.6,2.47){\circle{0.15}}
\put(7.8,2.4){\makebox(0,0)[c]{\scriptsize B}}
\put(7.6,0.5){\line(0,1){0.2}}
\put(7.6,0.25){\makebox(0,0)[c]{\scriptsize $\tau$}}
\put(6.85,1.35){\circle*{0.15}}
\put(6.85,0.5){\line(0,1){0.2}}
\put(6.85,0.25){\makebox(0,0)[c]{\scriptsize $\tau-\phi$}}
\multiput(6.85,1.35)(0,0.05){25}{\line(0,1){0.02}}
\multiput(6.85,2.6)(0.05,0){15}{\line(1,0){0.02}}
\put(7.55,2.6){\vector(1,0){0.01}}
\put(6.75,2){\vector(0,1){0.6}}
\put(6.75,2){\vector(0,-1){0.65}}
\put(6.65,2){\makebox(0,0)[c]{\scriptsize $\varepsilon$}}
\put(6.5,2.6){\line(1,0){0.15}}
\put(5.65,2.6){\makebox(0,0)[c]{\scriptsize$f(\tau-\phi)+\varepsilon$}}
\put(8.5,-0.5){\makebox(0,0)[c]{\footnotesize (e)}}

\put(12.5,0.5){\line(1,0){4}}
\put(12.5,4.5){\line(1,0){4}}
\put(12.5,0.5){\line(0,1){4}}
\put(16.5,0.5){\line(0,1){4}}
\qbezier(12.5,0.5)(13.5,3.5)(16.5,4.5)
\put(14.35,3.22){\circle*{0.15}}
\put(14.45,3){\makebox(0,0)[c]{\scriptsize A}}
\put(16.1,4.35){\circle{0.15}}
\put(16.1,4.1){\makebox(0,0)[c]{\scriptsize B}}
\put(14.35,0.5){\line(0,1){0.2}}
\put(14.35,0.25){\makebox(0,0)[c]{\scriptsize $\tau+\phi$}}
\multiput(14.35,3.45)(0,0.05){18}{\line(0,1){0.02}}
\multiput(14.35,4.35)(0.05,0){34}{\line(1,0){0.02}}
\put(16.05,4.35){\vector(1,0){0.01}}
\put(14.2,4){\vector(0,1){0.35}}
\put(14.2,4){\vector(0,-1){0.7}}
\put(14.1,3.8){\makebox(0,0)[c]{\scriptsize $\varepsilon$}}
\put(12.5,4.35){\line(1,0){0.15}}
\put(11.65,4.35){\makebox(0,0)[c]{\scriptsize$f(\tau+\phi)+\varepsilon$}}
\put(14.35,3.35){\circle{0.15}}
\put(14.5,-0.5){\makebox(0,0)[c]{\footnotesize (f)}}
\end{picture}

\vspace{3ex} {\footnotesize Fig.3: A system of two oscillators (A
and B) governed by $x=f(\varphi)$, with $\tau$ and $\varepsilon$
satisfying $\varepsilon<f(\tau)<1$. The initial phases are chosen
to be $\varphi_{A}(0)=1-\phi$ and $\varphi_{B}(0)=1$ such that
$0<\phi<\tau$ and $f(\tau-\phi)+\varepsilon=f(\tau)$. (a) The
state of the system at $t=0$. At the moment oscillator B emits a
spike. (b) The state of the system just after $t=0$. (c) The state
of the system at $t=\phi$. At the moment oscillator A emits a
spike. (d) The state of the system just after $t=\phi$. (e) The
state of the system at $t=\tau$. At the moment the spike emitted
by B at $t=0$ reaches A. (f) The state of the system at
$t=\tau+\phi$. At the moment the spike emitted by A at $t=\phi$
reaches B. Although equality $\varphi_{A}(s)=\varphi_{B}(s)$ holds
for all $s\in[\tau,\tau+\phi)$, (f) shows that the system has not
been synchronized at $t=\tau+\phi$.}

\begin{proposition}\quad 
$F_m(\theta)$ has the following properties:
\begin{list}
{{\upshape (\alph{mycounter})}\hfill}
{\setlength{\topsep}{0ex}
 \setlength{\parskip}{0ex}
 \setlength{\itemsep}{0.2ex}
 \setlength{\parsep}{0ex}
 \setlength{\leftmargin}{3.75ex}
 \setlength{\labelwidth}{3.5ex}
 \setlength{\labelsep}{0.25ex}
 \usecounter{mycounter}}
\item $0\leq F_m(\theta)\leq1$ for all $0\leq\theta\leq1$ and all
$m\in \mathbb{Z}^{+}$.

\item If $\theta_1\leq\theta_2$, then $F_m(\theta_1)\leq
F_m(\theta_2)$. In particular, if $\theta_1<\theta_2$ and
$F_m(\theta_1)<1$, then\\ $F_m(\theta_1)<F_m(\theta_2)$.
{\upshape(}the monotonicity with respect to $\theta${\upshape)}

\item If $m_1\leq m_2$, then $F_{m_1}(\theta)\leq
F_{m_2}(\theta)$. In particular, if $m_1<m_2$ and
$F_{m_1}(\theta)<1$, then $F_{m_1}(\theta)<F_{m_2}(\theta)$.
{\upshape(}the monotonicity with respect to $m${\upshape)}

\item When $0<F_m(\theta)<1$, the derivative of $F_m(\theta)$ with
respect to $\theta$ exists, and
$\frac{\mathrm{d}F_m(\theta)}{\mathrm{d}\theta}>1$.

\item If $\theta, \delta\geq0$ and $F_m(\theta+\delta)<1$,
then $F_m(\theta)+\delta\leq F_m(\theta+\delta)$\\
and the equal sign holds if and only if $\delta=0$ or $m=0$, i.e.,
$m\delta=0$.

\item $F_n(F_m(\theta))=F_{m+n}(\theta)$ for all
$0\leq\theta\leq1$ and all $m,n\in \mathbb{Z}^{+}$.

\item If $\theta,\delta_1,\ldots,\delta_n\geq0$ and
$F_{m_1+\cdots+m_n}(\theta+\delta_1+\cdots+\delta_n)<1$,
then\\
$F_{m_n}(\cdots
(F_{m_2}(F_{m_1}(\theta)+\delta_1)+\delta_2)+\cdots)+\delta_n
\leq F_{m_1+\cdots+m_n}(\theta+\delta_1+\cdots+\delta_n)$\\
and the equal sign holds if and only if
$m_1\delta_1=\cdots=m_n\delta_n=0$.
\end{list}
\end{proposition}

\begin{remark}\quad 
Property {\upshape (d)} comes from the concavity assumption of
function $f$. Property {\upshape (e)} follows from {\upshape (d)},
{\upshape (f)} is clear, and {\upshape (g)} follows from {\upshape
(e)} and {\upshape (f)}.
\end{remark}

\begin{proposition}\quad 
$D_i(t)$ and $D_i(t^{+})$ have the following properties:
\begin{list}
{{\upshape (\alph{mycounter})}\hfill} {\setlength{\topsep}{0ex}
 \setlength{\parskip}{0ex}
 \setlength{\itemsep}{0.2ex}
 \setlength{\parsep}{0ex}
 \setlength{\leftmargin}{4ex}
 \setlength{\labelwidth}{3.5ex}
 \setlength{\labelsep}{0.5ex}
 \usecounter{mycounter}}
\item For each $\eta\in D_i(t)$, we have $0\leq\eta<\tau$.

\item[{\upshape(a')}] For each $\eta\in D_i(t^{+})$, we have
$0<\eta\leq\tau$.

\item For any $\eta_1,\eta_2\in D_i(t)$, we have
$|\eta_1-\eta_2|<\tau$.

\item[{\upshape(b')}] For any $\eta_1,\eta_2\in D_i(t^{+})$, we
have $|\eta_1-\eta_2|<\tau$.

\item If $\eta\in D_i(t)$, then $t-[\tau-\eta]\in T_i$, i.e.,
$\varphi_i(t-[\tau-\eta])=1$.

\item[{\upshape(c')}] If $\eta\in D_i(t^{+})$, then
$t-[\tau-\eta]\in T_i$, i.e., $\varphi_i(t-[\tau-\eta])=1$.

\item If $\eta\in D_i(t)$ and $0<\eta<\tau$, then $\eta\in
D_i(t^{+})$ too.

\item[{\upshape(d')}] If $\eta\in D_i(t^{+})$ and $0<\eta<\tau$,
then $\eta\in D_i(t)$ too.

\item If $\eta\in D_i(t)$, then $\eta-\omega\in D_i(t+\omega)$,
where $\omega\in(\eta-\tau,\eta]$.

\item[{\upshape(e')}] If $\eta\in D_i(t^{+})$, then
$\eta-\omega\in D_i((t+\omega)^{+})$, where
$\omega\in[\eta-\tau,\eta)$.
\end{list}
\end{proposition}

\section{Properties of the pulse-coupled oscillator network
with delayed excitatory coupling under the assumptions (A1) and
(A2)}\quad It is difficult to prove the main theorem directly. The
following lemmas are needed. Lemma 1 shows that under the
assumptions (A1) and (A2), the time between any two firing
activities of any oscillator is longer than $2\tau$. This result is
reasonable to fireflies, because the transmission delays of
fireflies are relatively short and fireflies can not fire twice in
such a short period of time. By Lemma 1, we obtain Lemma 2, which
shows that under the assumptions (A1) and (A2), both $D_i(t)$ and
$D_i(t^{+})$ have at most one element. In Lemmas 3-5, three criteria
of synchronization of two oscillators are given. In Lemma 6 and
Lemma 7, two special cases are discussed, respectively. The two
lemmas show that the network can not achieve complete
synchronization even if it has been divided into two cliques.

\begin{lemma}\quad 
If an oscillator fires at time $t_1$ and $t_2$ with $t_1\neq t_2$,
then $|t_1-t_2|>2\tau$.
\end{lemma}
{\bf Proof:} Assume $t_2>t_1$ and that the oscillator is the $i$th
oscillator. Then we have $\varphi_i(t_1^{+})=0$ and
$\varphi_i(t_2)=1$. For $t_1,t_2$, there must be
$t_3,t_4\in[t_1,t_2]$ such that the $i$th oscillator fires at
$t_3,t_4$ and does not fire in $(t_3,t_4)$. If we can prove
$|t_3-t_4|>2\tau$, then $|t_1-t_2|>2\tau$ must hold. Therefore,
without loss of generality, we can assume that the $i$th oscillator
does not fire in time $(t_1,t_2)$.

Suppose $t_2-t_1\leq2\tau$.

{\bf Case 1}: $t_2<\tau$.

Because of the presence of the transmission delay and the assumption
(A1), no spikes can reach the $i$th oscillator in time
$(t_1,t_2]\subseteq(0,\tau)$. Then, according to Proposition 1,
$\varphi_i(t_2)=\varphi_i(t_2^{-})=\varphi_i(t_1^{+})+(t_2-t_1)=t_2-t_1<\tau<1$,
which contradicts that $\varphi_i(t_2)=1$. So, this case is
impossible.

{\bf Case 2}: $t_2=\tau$.

Because of the presence of the transmission delay and the assumption
(A1), only the spikes that the other oscillators emit at time $0$
can reach the $i$th oscillator in time $(t_1,t_2]$. Assume that
among the other oscillators, there are $m$ oscillators firing at
time $0$. (Here, $m$ can be equal to zero. It means that among the
other oscillators there are no oscillators firing at time 0. For the
case of $m=0$, the following proof holds true.) It is obvious that
$0\leq m\leq N-1$. According to Proposition 1 and the assumption
(A2), we have
\begin{eqnarray*}
\varphi_i(t_2^{-})=\varphi_i(t_1^{+})+(t_2-t_1)=\tau-t_1\leq\tau
\end{eqnarray*}
and
\begin{eqnarray*}
\varphi_i(t_2) &  =  & f^{-1}(\min[1,f(\varphi_i(t_2^{-}))+m\varepsilon])\\
               &  =  & f^{-1}(\min[1,f(t_2-t_1)+m\varepsilon])\\
               & \leq& f^{-1}(\min[1,f(\tau)+(N-1)\varepsilon])\\
               & \leq& f^{-1}(\min[1,f(2\tau)+N\varepsilon])\\
               &  =  & f^{-1}(f(2\tau)+N\varepsilon)\\
               &  <  & f^{-1}(1)=1
\end{eqnarray*}
which contradicts that $\varphi_i(t_2)=1$. So, this case is
impossible.

{\bf Case 3}: $t_2>\tau$.

Assume that there are $m=m_1+m_2+\cdots+m_l$ spikes which reach
the $i$th oscillator in time $(t_1,t_2]$:
\begin{eqnarray*}
&&m_1 \mbox{ spikes reach at time } t_1+\delta_1;\\
&&m_2 \mbox{ spikes reach at time } t_1+\delta_1+\delta_2;\\
&&\cdots\cdots\\
&&m_l \mbox{ spikes reach at time } t_1+\delta_1+\cdots+\delta_l.
\end{eqnarray*}
where $\delta_1,\ldots,\delta_l>0$, and
$\delta_1+\cdots+\delta_l\leq t_2-t_1$. (Here, $m$ can be equal to
zero. It means that there are no spikes reaching the $i$th
oscillator in time $(t_1,t_2]$. When $m=0$, we can let $l=1$,
$m_l=0$ and $\delta_l$ be any value in the range $(0,t_2-t_1]$.
Hence the following proof holds true for the case of $m=0$.)

Because the $i$th oscillator does not fire in time $(t_1,t_2)$,
according to Proposition 1 we have
{\footnotesize
\begin{eqnarray*}
&&  \varphi_i([t_1+\delta_1]^{-})=\varphi_i(t_1^{+})+\delta_1=\delta_1\\
&&  \varphi_i(t_1+\delta_1)=F_{m_1}(\varphi_i([t_1+\delta_1]^{-}))=F_{m_1}(\delta_1)\\
&&
\varphi_i([t_1+\delta_1+\delta_2]^{-})=\varphi_i(t_1+\delta_1)+\delta_2
    =F_{m_1}(\delta_1)+\delta_2\\
&&
\varphi_i(t_1+\delta_1+\delta_2)=F_{m_2}(\varphi_i([t_1+\delta_1+\delta_2]^{-}))
    =F_{m_2}(F_{m_1}(\delta_1)+\delta_2)\\
&&  \cdots \cdots\\
&&
\varphi_i([t_1+\delta_1+\cdots+\delta_l]^{-})=\varphi_i(t_1+\delta_1+\cdots+\delta_{l-1})+\delta_l
    =F_{m_{l-1}}(\cdots F_{m_2}(F_{m_1}(\delta_1)+\delta_2)\cdots+\delta_{l-1})+\delta_l\\
&&
\varphi_i(t_1+\delta_1+\cdots+\delta_l)=F_{m_l}(\varphi_i([t_1+\delta_1+\cdots+\delta_l]^{-}))
    =F_{m_l}(F_{m_{l-1}}(\cdots F_{m_2}(F_{m_1}(\delta_1)+\delta_2)\cdots+\delta_{l-1})+\delta_l)
\end{eqnarray*}}

\vspace{-5ex}

Because of the presence of the transmission delay and the
assumption (A1), the spikes which reach the $i$th oscillator in
time $(t_1,t_2]$ are the ones that the other oscillators emit in
time $(t_1-\tau,t_2-\tau]\cap[0,+\infty)$. If all the other
oscillators fire not more than once in
$(t_1-\tau,t_2-\tau]\cap[0,+\infty)$, i.e., $m\leq N-1$, then it
follows from the assumption (A2) and Proposition 2(g) that

\vspace{-4ex} {\footnotesize
\begin{eqnarray*}
&&  \varphi_i(t_2)=\varphi_i(t_1+\delta_1+\cdots+\delta_l)+[t_2-(t_1+\delta_1+\cdots+\delta_l)]\\
&&  \hspace{29pt}=F_{m_l}(F_{m_{l-1}}(\cdots F_{m_2}(F_{m_1}(\delta_1)+\delta_2)\cdots+\delta_{l-1})+\delta_l)
                  +[t_2-(t_1+\delta_1+\cdots+\delta_l)]\\
&&  \hspace{29pt}\leq F_{m_1+\cdots+m_l}(\delta_1+\cdots+\delta_l+[t_2-(t_1+\delta_1+\cdots+\delta_l)])\\
&&  \hspace{29pt}=F_{m_1+\cdots+m_l}(t_2-t_1)\\
&&  \hspace{29pt}\leq F_N(2\tau)<1
\end{eqnarray*}}
which contradicts that $\varphi_i(t_2)=1$. Hence there must be some
oscillator firing more than once in
$(t_1-\tau,t_2-\tau]\cap[0,+\infty)$. It means that there must be
time $t_3,t_4$ with $t_3<t_4$ and
$t_3,t_4\in(t_1-\tau,t_2-\tau]\cap[0,+\infty)$ such that some
oscillator fires at $t_3,t_4$ and does not fire in $(t_3,t_4)$. If
$t_4>\tau$, as above, there must be time $t_5,t_6$ with $t_5<t_6$
and $t_5,t_6\in(t_3-\tau,t_4-\tau]\cap[0,+\infty)$ such that some
oscillator fires at $t_5,t_6$ and does not fire in $(t_5,t_6)$. If
$t_6>\tau$, $\cdots$ This leads to a finite sequence:
\begin{eqnarray*}
&(t_1,t_2)\rightarrow(t_3,t_4)\rightarrow\cdots\rightarrow(t_{2k-1},t_{2k})
\rightarrow\cdots\rightarrow(t_{2n-1},t_{2n})
\end{eqnarray*}
which satisfies that
\begin{eqnarray*}
&& t_{2(k-1)-1}-\tau<t_{2k-1}<t_{2k}\leq t_{2(k-1)}-\tau, k=2,\ldots,n,\\
&& t_{2k-1}\geq0, k=1,\ldots,n,\\
&& t_{2k}>\tau, k=1,\ldots,n-1 \mbox{, and } t_{2n}\leq\tau.
\end{eqnarray*}
The last term of the sequence means that some oscillator fires at
$t_{2n-1}$ and $t_{2n}$ and does not fire in $(t_{2n-1},t_{2n})$,
with $0\leq t_{2n-1}<t_{2n}\leq\tau$. By case 1 and case 2, this
is impossible. So, case 3 is impossible, too.

Therefore, we have $t_2-t_1>2\tau$, which completes the proof.
\hfill $\square$

By Lemma 1, we get Lemma 2.

\begin{lemma}\quad 
Under the assumptions {\upshape(A1)} and {\upshape(A2)}, $D_i(t)$
and $D_i(t^{+})$ have the following properties:
\begin{list}
{{\upshape (\alph{mycounter})}\hfill} {\setlength{\topsep}{0ex}
 \setlength{\parskip}{0ex}
 \setlength{\itemsep}{0.2ex}
 \setlength{\parsep}{0ex}
 \setlength{\leftmargin}{3.75ex}
 \setlength{\labelwidth}{3.5ex}
 \setlength{\labelsep}{0.25ex}
 \usecounter{mycounter}}
\item $D_i(t)=\{\eta\}$ with $0\leq\eta<\tau$ or
$D_i(t)=\emptyset$, for all $1\leq i\leq N$ and all $t\geq0$;
i.e., set $D_i(t)$ has at most one element for all $1\leq i\leq N$
and all $t\geq0$.

\item $D_i(t^{+})=\{\eta\}$ with $0<\eta\leq\tau$ or
$D_i(t^{+})=\emptyset$, for all $1\leq i\leq N$ and all $t\geq0$;
i.e., set $D_i(t^{+})$ has at most one element for all $1\leq
i\leq N$ and all $t\geq0$.

\item If $D_i(t)\neq\emptyset$, then $0<\varphi_i(t)<1$.

\item If $D_i(t)=\{\eta_1\}$, $D_j(t)=\{\eta_2\}$ and
$\eta_1<\eta_2$, then $\varphi_i(t^{-})>\varphi_j(t^{-})$.
\end{list}
\end{lemma}
{\bf Proof:} (a) Assume that there are $\eta_1,\eta_2\in D_i(t)$
with $\eta_1\neq\eta_2$. Then by Proposition 3(c) we have
$\varphi_i(t-[\tau-\eta_1])=\varphi_i(t-[\tau-\eta_2])=1$, i.e.,
the $i$th oscillator fires at time $t-[\tau-\eta_1]$ and
$t-[\tau-\eta_2]$. However, it follows from Proposition 3(b) that
$\big|[t-(\tau-\eta_1)]-[t-(\tau-\eta_2)]\big|=|\eta_1-\eta_2|<\tau$.
This contradicts Lemma 1. So, for any $\eta_1,\eta_2\in D_i(t)$,
the equality $\eta_1=\eta_2$ must hold. It shows that
$D_i(t)=\{\eta\}$ with $0\leq\eta<\tau$, or $D_i(t)=\emptyset$.

(b) The proof is similar to that of (a).

(c) If $D_i(t)\neq\emptyset$, there must be a $\eta\in D_i(t)$
with $0\leq\eta<\tau$. By Proposition 3(c), we have
$\varphi_i(t-[\tau-\eta])=1$. Since
$t-[t-(\tau-\eta)]=\tau-\eta<2\tau$, by Lemma 1, the $i$th
oscillator never fires at time $t$. It means $\varphi_i(t)<1$.
Thus, from Remark 1, we have $0<\varphi_i(t)<1$.

(d) Assume that there are $m=m_1+m_2+\cdots+m_l$ spikes which
reach the $i$th oscillator in time $(t-(\tau-\eta_2),t)$:
\begin{eqnarray*}
&&m_1 \mbox{ spikes reach at time }t-(\tau-\eta_2)+\delta_1;\\
&&m_2 \mbox{ spikes reach at time }t-(\tau-\eta_2)+\delta_1+\delta_2;\\
&&\cdots\cdots\\
&&m_l \mbox{ spikes reach at time }t-(\tau-\eta_2)+\delta_1+\cdots+\delta_l.
\end{eqnarray*}
where $\delta_1,\ldots,\delta_l>0$, and
$\delta_1+\cdots+\delta_l<\tau-\eta_2$. (Here, $m$ can be equal to
zero. It means that there are no spikes reaching the $i$th
oscillator in time $(t-(\tau-\eta_2),t)$. When $m=0$, we can let
$l=1$, $m_l=0$ and $\delta_l$ be any value in the range
$(0,\tau-\eta_2)$. Hence the following proof holds true for the case
of $m=0$.)

Since $D_i(t)=\{\eta_1\}$, $D_j(t)=\{\eta_2\}$ and $\eta_1<\eta_2$,
by Proposition 3(c) and 3(e) we have that
$\varphi_j(t-(\tau-\eta_2))=1$, $D_i(s)=\{\eta_1+(t-s)\}$ for all
$s\in[t-(\tau-\eta_2),t]$ and $D_j(s)=\{\eta_2+(t-s)\}$ for all
$s\in(t-(\tau-\eta_2),t]$.

From $D_i(s)=\{\eta_1+(t-s)\}\neq\{0\}$ and
$D_j(s)=\{\eta_2+(t-s)\}\neq\{0\}$ for all
$s\in(t-(\tau-\eta_2),t)$, it follows that oscillator $i$ and $j$
can not receive the spikes of each other in time
$(t-(\tau-\eta_2),t)$, which means that the spikes reaching the
$j$th oscillator in time $(t-(\tau-\eta_2),t)$ are the $m$ spikes
above.

Because $D_i(s)=\{\eta_1+(t-s)\}\neq\emptyset$ and
$D_j(s)=\{\eta_2+(t-s)\}\neq\emptyset$ for all
$s\in(t-(\tau-\eta_2),t)$, from Lemma 2(c), it follows that
$0<\varphi_i(s)<1$ and $0<\varphi_j(s)<1$ for all
$s\in(t-(\tau-\eta_2),t)$, that is, both oscillator $i$ and $j$ do
not fire in time $(t-(\tau-\eta_2),t)$.

Since $D_i(t-(\tau-\eta_2))=\{\tau+\eta_1-\eta_2\}\neq\emptyset$, by
Lemma 2(c) we have $0<\varphi_i(t-(\tau-\eta_2))<1$. Denote
$\phi=\varphi_i(t-(\tau-\eta_2))$.

Therefore, according to Proposition 1, we have {\footnotesize
\begin{eqnarray*}
&&\varphi_i([t-(\tau-\eta_2)+\delta_1]^{-})=\phi+\delta_1\\
&&\varphi_i(t-(\tau-\eta_2)+\delta_1)=F_{m_1}(\phi+\delta_1)\\
&&\varphi_i([t-(\tau-\eta_2)+\delta_1+\delta_2]^{-})=F_{m_1}(\phi+\delta_1)+\delta_2\\
&&\varphi_i(t-(\tau-\eta_2)+\delta_1+\delta_2)=F_{m_2}(F_{m_1}(\phi+\delta_1)+\delta_2)\\
&&\cdots\cdots\\
&&\varphi_i([t-(\tau-\eta_2)+\delta_1+\cdots+\delta_l]^{-})=
  F_{m_{l-1}}(\cdots F_{m_2}(F_{m_1}(\phi+\delta_1)+\delta_2)\cdots+\delta_{l-1})+\delta_l\\
&&\varphi_i(t-(\tau-\eta_2)+\delta_1+\cdots+\delta_l)=
  F_{m_l}(F_{m_{l-1}}(\cdots F_{m_2}(F_{m_1}(\phi+\delta_1)+\delta_2)\cdots+\delta_{l-1})+\delta_l)\\
\end{eqnarray*}}
and
{\footnotesize
\begin{eqnarray*}
&&\varphi_j([t-(\tau-\eta_2)+\delta_1]^{-})=\delta_1\\
&&\varphi_j(t-(\tau-\eta_2)+\delta_1)=F_{m_1}(\delta_1)\\
&&\varphi_j([t-(\tau-\eta_2)+\delta_1+\delta_2]^{-})=F_{m_1}(\delta_1)+\delta_2\\
&&\varphi_j(t-(\tau-\eta_2)+\delta_1+\delta_2)=F_{m_2}(F_{m_1}(\delta_1)+\delta_2)\\
&&\cdots\cdots\\
&&\varphi_j([t-(\tau-\eta_2)+\delta_1+\cdots+\delta_l]^{-})=
  F_{m_{l-1}}(\cdots F_{m_2}(F_{m_1}(\delta_1)+\delta_2)\cdots+\delta_{l-1})+\delta_l\\
&&\varphi_j(t-(\tau-\eta_2)+\delta_1+\cdots+\delta_l)=
  F_{m_l}(F_{m_{l-1}}(\cdots F_{m_2}(F_{m_1}(\delta_1)+\delta_2)\cdots+\delta_{l-1})+\delta_l)\\
\end{eqnarray*}}
By Proposition 2(b), we have {\footnotesize
\begin{eqnarray*}
&&\varphi_i(t-(\tau-\eta_2)+\delta_1)=F_{m_1}(\phi+\delta_1)>F_{m_1}(\delta_1)=\varphi_j(t-(\tau-\eta_2)+\delta_1)\\
&&\varphi_i(t-(\tau-\eta_2)+\delta_1+\delta_2)=F_{m_2}(F_{m_1}(\phi+\delta_1)+\delta_2)
  >F_{m_2}(F_{m_1}(\delta_1)+\delta_2)=\varphi_j(t-(\tau-\eta_2)+\delta_1+\delta_2)\\
&&\cdots\cdots\\
&&\varphi_i(t-(\tau-\eta_2)+\delta_1+\cdots+\delta_l)
  =F_{m_l}(\cdots F_{m_2}(F_{m_1}(\phi+\delta_1)+\delta_2)\cdots+\delta_l)\\
&&\hspace{137pt}>F_{m_l}(\cdots F_{m_2}(F_{m_1}(\delta_1)+\delta_2)\cdots+\delta_l)
  =\varphi_j(t-(\tau-\eta_2)+\delta_1+\cdots+\delta_l)
\end{eqnarray*}}
Therefore,
\begin{eqnarray*}
\varphi_i(t^{-})
& =  & \varphi_i(t-(\tau-\eta_2)+\delta_1+\cdots+\delta_l)
       +\big\{t-[t-(\tau-\eta_2)+\delta_1+\cdots+\delta_l]\big\}\\
& >  & \varphi_j(t-(\tau-\eta_2)+\delta_1+\cdots+\delta_l)
       +\big\{t-[t-(\tau-\eta_2)+\delta_1+\cdots+\delta_l]\big\}\\
& =  & \varphi_j(t^{-})\mbox{\hspace{32em}$\square$}
\end{eqnarray*}

In the following (Lemmas 3-5), three synchronization criteria of two
oscillators will be given.

\begin{lemma}\quad 
Oscillator $i$ and $j$ have been synchronized at time $t_0$, if
and only if the two oscillators satisfy
$\varphi_i(t_0)=\varphi_j(t_0)$ and $D_i(t_0^{+})=D_j(t_0^{+})$.
\end{lemma}
{\bf Proof:} (1) Suppose that $\varphi_i(t_0)=\varphi_j(t_0)$ and
$D_i(t_0^{+})=D_j(t_0^{+})$.

It means that the two oscillators have acted as one, because their
dynamics are identical and they are coupled in the same way to all
the other oscillators. That is, oscillator $i$ and $j$ have been
synchronized at $t_0$.

(2) Suppose that oscillator $i$ and $j$ have been synchronized at
$t_0$, namely, $\varphi_i(t)=\varphi_j(t)$ for all $t\geq t_0$.

We only need to prove $D_i(t_0^{+})=D_j(t_0^{+})$. Assume
$D_i(t_0^{+})\neq D_j(t_0^{+})$.

It is clear that $D_i(t_0^{+})\neq\emptyset$ or
$D_j(t_0^{+})\neq\emptyset$. Without loss of generality, we let
$D_i(t_0^{+})\neq\emptyset$. By Lemma 2(b), we have
$D_i(t_0^{+})=\{\eta\}\neq D_j(t_0^{+})$ with $0<\eta\leq\tau$.
Then by Proposition 3, we get $D_i(t_0+\eta)=\{0\}\neq
D_j(t_0+\eta)$, which means that a spike of oscillator $i$ will
reach oscillator $j$ at time $t_0+\eta$, but no spikes of
oscillator $j$ will reach oscillator $i$ at time $t_0+\eta$. We
let $m_i$ be the number of the spikes which reach oscillator $i$
at time $t_0+\eta$, and $m_j$ be the number of the spikes which
reach oscillator $j$ at time $t_0+\eta$. So, the equality
$m_i=m_j-1$ holds, because the two oscillators are coupled in the
same way to all the other oscillators. It is obvious that
$\varphi_i([t_0+\eta]^{-})=\varphi_j([t_0+\eta]^{-})$, because
$\varphi_i(t)=\varphi_j(t)$ for all $t\geq t_0$. Then, according
to Proposition 1, we have
$\varphi_i(t_0+\eta)=f^{-1}(\min[1,f(\varphi_i([t_0+\eta]^{-}))+m_i\varepsilon])
\neq
f^{-1}(\min[1,f(\varphi_j([t_0+\eta]^{-}))+m_j\varepsilon])=\varphi_j(t_0+\eta)$,
unless $f(\varphi_i([t_0+\eta]^{-}))+m_i\varepsilon\geq1$ and
$f(\varphi_j([t_0+\eta]^{-}))+m_j\varepsilon\geq1$, i.e.,
$\varphi_i(t_0+\eta)=\varphi_j(t_0+\eta)=1$. However, by
$D_i(t_0+\eta)=\{0\}\neq\emptyset$ and Lemma 2(c), we have
$\varphi_i(t_0+\eta)<1$. This is a contradiction. Thus,
$D_i(t_0^{+})=D_j(t_0^{+})$. \hfill$\square$

\vspace{-1.5ex}
\begin{lemma}\quad 
If $D_i(t_0)=D_j(t_0)=\{\eta\}$ with $0\leq\eta<\tau$, then
oscillator $i$ and $j$ have been synchronized at time
$t_0-(\tau-\eta)$; conversely, if there exists a $t_1<t_0$ such that
oscillator $i$ and $j$ have been synchronized at time $t_1$, then
$D_i(t_0)=D_j(t_0)$.
\end{lemma}
\vspace{-1.5ex}{\bf Proof:} (1) Suppose that
$D_i(t_0)=D_j(t_0)=\{\eta\}$ with $0\leq\eta<\tau$.

From Proposition 3, we have
$D_i([t_0-(\tau-\eta)]^{+})=D_j([t_0-(\tau-\eta)]^{+})=\{\tau\}$ and
$\varphi_i(t_0-(\tau-\eta))=\varphi_j(t_0-(\tau-\eta))=1$. By Lemma
3, oscillator $i$ and $j$ have been synchronized at time
$t_0-(\tau-\eta)$.

(2) Suppose that there exists a $t_1<t_0$ such that oscillator $i$
and $j$ have been synchronized at time $t_1$.

Assume $D_i(t_0)\neq D_j(t_0)$. It is clear that
$D_i(t_0)\neq\emptyset$ or $D_j(t_0)\neq\emptyset$. Without loss of
generality, we let $D_i(t_0)\neq\emptyset$. By Lemma 2(a), we have
$D_i(t_0)=\{\eta\}\neq D_j(t_0)$ with $0\leq \eta<\tau$. If follows
from Proposition 3 that $D_i(t_0+\eta)=\{0\}\neq D_j(t_0+\eta)$.
Since oscillator $i$ and $j$ have been synchronized at time $t_1$
with $t_1<t_0$, by an argument similar to that used in the second
part of the proof of Lemma 3, we can obtain
$\varphi_i(t_0+\eta)=\varphi_j(t_0+\eta)=1$. However, by
$D_i(t_0+\eta)=\{0\}\neq\emptyset$ and Lemma 2(c), we have
$\varphi_i(t_0+\eta)<1$. This is a contradiction. Thus,
$D_i(t_0)=D_j(t_0)$. \hfill$\square$

\vspace{-1.5ex}
\begin{lemma}\quad 
Oscillator $i$ and $j$ have been synchronized at time $t_0>0$ and
they have not been synchronized at time $t$ for all $0\leq t<t_0$,
if and only if  the two oscillators satisfy one of the following
cases:
\begin{list}
{{\upshape (\alph{mycounter})}\hfill}
{\setlength{\topsep}{0ex}
 \setlength{\parskip}{0ex}
 \setlength{\itemsep}{0.2ex}
 \setlength{\parsep}{0ex}
 \setlength{\leftmargin}{3.75ex}
 \setlength{\labelwidth}{3.5ex}
 \setlength{\labelsep}{0.25ex}
 \usecounter{mycounter}}
\item $\varphi_i(t_0^{-})\neq\varphi_j(t_0^{-})$,
$D_i(t_0)=D_j(t_0)=\emptyset$ and
$\varphi_i(t_0)=\varphi_j(t_0)=1$.

\item $\varphi_i(t_0^{-})\neq\varphi_j(t_0^{-})$,
$D_i(t_0)=\{0\}$, $D_j(t_0)=\emptyset$ and
$\varphi_i(t_0)=\varphi_j(t_0)<1$.

\item $\varphi_i(t_0^{-})\neq\varphi_j(t_0^{-})$,
$D_i(t_0)=\emptyset$, $D_j(t_0)=\{0\}$ and
$\varphi_i(t_0)=\varphi_j(t_0)<1$.
\end{list}
\end{lemma}
\vspace{-1.5ex}{\bf Proof:}  (1) Suppose that oscillator $i$ and $j$
have been synchronized at $t_0$ and that they have not been
synchronized at $t$ for all $0\leq t<t_0$.

We have $\varphi_i(t_0)=\varphi_j(t_0)$, because oscillator $i$ and
$j$ have been synchronized at $t_0$. Assume that
$\varphi_i(t_0^{-})=\varphi_j(t_0^{-})$. From Lemma 1, there must be
a $\delta>0$ such that no oscillators fire and no spikes reach in
time $(t_0-\delta,t_0)$. It means that
$\mathrm{d}\varphi_i(t)/\mathrm{d}t=\mathrm{d}\varphi_j(t)/\mathrm{d}t=1$
for all $t\in(t_0-\delta,t_0)$. Thus, we get
$\varphi_i(t)=\varphi_j(t)$ for all $t\geq t_0-\delta/2$. This
contradicts that oscillator $i$ and $j$ have not been synchronized
at $t$ for all $0\leq t<t_0$. So,
$\varphi_i(t_0^{-})\neq\varphi_j(t_0^{-})$.

If no spikes reach both oscillator $i$ and $j$ at time $t_0$, then
it follows from $\varphi_i(t_0^{-})\neq\varphi_j(t_0^{-})$ that
$\varphi_i(t_0)=\varphi_i(t_0^{-})\neq\varphi_j(t_0^{-})=\varphi_j(t_0)$.
It contradicts that oscillator $i$ and $j$ have been synchronized at
$t_0$. Therefore, there must be a spike which reaches at least one
of the two oscillators at time $t_0$, i.e, there must be an
oscillator $k$ such that $D_k(t_0)=\{0\}$. (This is an important
property of pulse-coupled networks.) Let $K$ be a set of
oscillators: oscillator $k'\in K$ if and only if
$D_{k'}(t_0)=\{0\}$. And let $m$ be the number of the oscillators in
$K$.

{\bf Case 1}: oscillator $i,j\not\in K$.

Then
$f^{-1}(\min[1,f(\varphi_i(t_0^{-}))+m\varepsilon])=\varphi_i(t_0)=
\varphi_j(t_0)=f^{-1}(\min[1,f(\varphi_j(t_0^{-}))+m\varepsilon])$.
But $\varphi_i(t_0^{-})\neq\varphi_j(t_0^{-})$. So,
$\varphi_i(t_0)=\varphi_j(t_0)=1$. By Lemma 2(c), we get
$D_i(t_0)=D_j(t_0)=\emptyset$. Then case (a) holds.

{\bf Case 2}: oscillator $i\in K$ and oscillator $j\not\in K$.

Then we have $D_i(t_0)=\{0\}\neq D_j(t_0)$. $D_i(t_0)=\{0\}$ shows
that $D_i(t_0^{+})=\emptyset$ and $\varphi_i(t_0)<1$. Since
oscillator $i$ and $j$ have been synchronized at $t_0$, it follows
from Lemma 3 that $D_j(t_0^{+})=D_i(t_0^{+})=\emptyset$, which
means $D_j(t_0)=\emptyset\,\,or\,\,\{0\}$. But
$D_j(t_0)\neq\{0\}$. So, $D_j(t_0)=\emptyset$. Then case (b)
holds.

{\bf Case 3}: oscillator $i\not\in K$ and oscillator $j\in K$.

It is similar to case 2. And case (c) holds.

{\bf Case 4}: oscillator $i,j\in K$.

Then $D_i(t_0)=D_j(t_0)=\{0\}$. By Lemma 4, oscillator $i$ and $j$
have been synchronized at time $t_0-\tau<t_0$, which contradicts
that oscillator $i$ and $j$ have not been synchronized at $t$ for
all $0\leq t<t_0$. Therefore, this case is impossible.

(2) Suppose that one of the three cases (a), (b) and (c) holds. Then
we have $\varphi_i(t_0)=\varphi_j(t_0)$ and
$D_i(t_0^{+})=D_j(t_0^{+})$. From Lemma 3, oscillator $i$ and $j$
have been synchronized at $t_0$. If there exists a $t'<t_0$ such
that oscillator $i$ and $j$ have been synchronized at $t'$, then
$\varphi_i(t_0^{-})=\varphi_j(t_0^{-})$. This is a contradiction.
So, oscillator $i$ and $j$ have not been synchronized at $t$ for all
$0\leq t<t_0$.$\square$

Two special cases are discussed in the following two lemmas, which
show that the network can not achieve synchronization even if it has
been divided into two cliques.

\begin{lemma}\quad 
If there exists a $t_0\geq 0$ such that
$\varphi_1(t_0)=\cdots=\varphi_m(t_0)=1-\varphi$
{\upshape($0<\varphi<1$)}, $D_1(t_0)=\cdots=D_m(t_0)=\emptyset$ and
$\varphi_{m+1}(t_0)=\cdots=\varphi_N(t_0)=1$,
$D_{m+1}(t_0)=\cdots=D_N(t_0)=\emptyset$, where $1\leq m<N$, then
the pulse-coupled oscillator network can not achieve
synchronization.
\end{lemma}
{\bf Proof:} Let $n=N-m$.

$\varphi_1(t_0)=\cdots=\varphi_m(t_0)=1-\varphi$
{\upshape($0<\varphi<1$)}, $D_1(t_0)=\cdots=D_m(t_0)=\emptyset$
and $\varphi_{m+1}(t_0)=\cdots=\varphi_N(t_0)=1$,
$D_{m+1}(t_0)=\cdots=D_N(t_0)=\emptyset$ show that oscillator
$1,\ldots,m$ have been synchronized at $t_0$, and oscillator
$m+1,\ldots,N$ have also been synchronized at $t_0$. Therefore, we
only need to discuss two oscillators, the $i$th oscillator and the
$j$th oscillator ($1\leq i\leq m$, $m+1\leq j\leq N$). The
synchronization of the whole network is equivalent to the
synchronization of oscillator $i$ and $j$.

According to the model rules, we get the following tables of
$\varphi_i(t)$, $D_i(t)$, $\varphi_j(t)$ and $D_j(t)$:

{\bf Case 1}: $\varphi<\tau$.
{\footnotesize
\begin{center}
\begin{tabular}{|l|l|l|l|l|}\hline
time $t$ & $\varphi_i(t)$ & $D_i(t)$ & $\varphi_j(t)$ &
$D_j(t)$\\
\hline \hline
$t_0$ & $1-\varphi$ & $\emptyset$ & 1 & $\emptyset$\\
\hline
$t_0^{+}$ & $1-\varphi$ & $\emptyset$ & 0 & $\{\tau\}$\\
\hline
$t_0+\varphi$ & 1 & $\emptyset$ & $\varphi$ & $\{\tau-\varphi\}$\\
\hline
$(t_0+\varphi)^{+}$ & 0 & $\{\tau\}$ & $\varphi$ & $\{\tau-\varphi\}$\\
\hline
$t_0+\tau$ & $F_n(\tau-\varphi)$ & $\{\varphi\}$ & $F_{n-1}(\tau)$ & $\{0\}$\\
\hline
$(t_0+\tau)^{+}$ & $F_n(\tau-\varphi)$ & $\{\varphi\}$ & $F_{n-1}(\tau)$ & $\emptyset$\\
\hline
$t_0+\tau+\varphi$ & $F_{m-1}(F_n(\tau-\varphi)+\varphi)$ &
$\{0\}$ &
$F_m(F_{n-1}(\tau)+\varphi)$ & $\emptyset$\\
\hline
$(t_0+\tau+\varphi)^{+}$ &
$F_{m-1}(F_n(\tau-\varphi)+\varphi)$ & $\emptyset$ &
$F_m(F_{n-1}(\tau)+\varphi)$ & $\emptyset$\\
\hline
$t_0+\tau+\varphi+$ &
$1-[F_m(F_{n-1}(\tau)+\varphi)$ & $\emptyset$ & 1 & $\emptyset$\\
$1-F_m(F_{n-1}(\tau)+\varphi)$ &
$-F_{m-1}(F_n(\tau-\varphi)+\varphi)]$ &&&\\
\hline
\end{tabular}
\end{center}}

\vspace{2ex} {\bf Case 2}: $\varphi\geq\tau$. {\footnotesize
\begin{center}
\begin{tabular}{|l|l|l|l|l|}\hline
time $t$ & $\varphi_i(t)$ & $D_i(t)$ & $\varphi_j(t)$ &
$D_j(t)$\\
\hline \hline
$t_0$ & $1-\varphi$ & $\emptyset$ & 1 & $\emptyset$\\
\hline
$t_0^{+}$ & $1-\varphi$ & $\emptyset$ & 0 & $\{\tau\}$\\
\hline
$t_0+\tau$ & $F_n(1-\varphi+\tau)$ &
$\emptyset$ & $F_{n-1}(\tau)$ & $\{0\}$\\
\hline
$(t_0+\tau)^{+}$ & $F_n(1-\varphi+\tau)$ &
$\emptyset$ & $F_{n-1}(\tau)$ & $\emptyset$\\
\hline
$t_0+\tau+$ & 1 & $\emptyset$ &
$1-[F_n(1-\varphi+\tau)$ & $\emptyset$\\
$1-F_n(1-\varphi+\tau)$ &&& $-F_{n-1}(\tau)]$ &\\
\hline
\end{tabular}
\end{center}}

\bigskip

Let
\begin{eqnarray*}
G_1(\theta)&=&F_m(F_{n-1}(\tau)+\theta)-F_{m-1}(F_n(\tau-\theta)+\theta)\qquad0\leq\theta<\tau\\
G_2(\theta)&=&F_n(F_{m-1}(\tau)+\theta)-F_{n-1}(F_m(\tau-\theta)+\theta)\qquad0\leq\theta<\tau\\
G_3(\theta)&=&F_n(1-\theta+\tau)-F_{n-1}(\tau)\qquad\qquad\qquad\qquad\,\,\,\,\,\tau\leq\theta<1\\
G_4(\theta)&=&F_m(1-\theta+\tau)-F_{m-1}(\tau)\qquad\qquad\qquad\qquad\,\,\,\tau\leq\theta<1
\end{eqnarray*}

Define map of the form
\begin{eqnarray*}
G(\theta,P,Q)=\left\{
  \begin{array}{cccc}
  (G_1(\theta),P,Q) & \qquad if\quad 0\leq\theta<\tau,P=m,Q=n\\
  (G_2(\theta),P,Q) & \qquad if\quad 0\leq\theta<\tau,P=n,Q=m\\
  (G_3(\theta),Q,P) & \qquad if\quad \tau\leq\theta<1,P=m,Q=n\\
  (G_4(\theta),Q,P) & \qquad if\quad \tau\leq\theta<1,P=n,Q=m
\end{array}\right.
\end{eqnarray*}
where $0\leq\theta<1$, $P=m \,\, or \,\, n$, $Q=m \,\, or \,\, n$.

Then the difference between the phase variables $\varphi_i$ and
$\varphi_j$ can be shown by the following sequence:
\begin{eqnarray*}
(\varphi,m,n) \rightarrow G(\varphi,m,n) \rightarrow
G^2(\varphi,m,n) \rightarrow \cdots
\end{eqnarray*}
If the pulse-coupled oscillator network can achieve synchronization,
then there must be a $p\in \mathbb{N}$ such that
$G^p(\varphi,m,n)=(0,m,n)\,\,or\,\,(0,n,m)$. In the following, we
will prove that such $p$ does not exist.

When $0<\theta<\tau$, we have $0<F_n(\tau-\theta)<1$,
$0<F_m(F_{n-1}(\tau)+\theta)<1$ and
$0<F_{m-1}(F_n(\tau-\theta)+\theta)<1$. From Proposition 2(d), it
follows that $F'_n(\tau-\theta)>1$, $F'_m(F_{n-1}(\tau)+\theta)>1$
and $F'_{m-1}(F_n(\tau-\theta)+\theta)>1$, for all $0<\theta<\tau$.
Therefore,
\begin{eqnarray*}
\frac{dG_1(\theta)}{d\theta}
=F'_m(F_{n-1}(\tau)+\theta)-F'_{m-1}(F_n(\tau-\theta)+\theta)\cdot[-F'_n(\tau-\theta)+1]>1
\end{eqnarray*}
for all $0<\theta<\tau$. Then from $ G_1(0)=0$, it follows that
$G_1(\theta)>0$ for all $0<\theta<\tau$.

When $\tau\leq\theta<1$, since $1-\theta+\tau>\tau$ and
$F_{n-1}(\tau)<1$, by Proposition 2(b) and 2(c) we have that
\begin{eqnarray*}
G_3(\theta)=F_n(1-\theta+\tau)-F_{n-1}(\tau)>0
\end{eqnarray*}
for all $\tau\leq\theta<1$.

Similarly, we can get that $G_2(\theta)>0$ for all $0<\theta<\tau$
and $G_4(\theta)>0$ for all $\tau\leq\theta<1$.

Thus, there does not exist $p\in \mathbb{N}$ such that
$G^p(\varphi,m,n)=(0,m,n)\,\,or\,\,(0,n,m)$. It means that the
pulse-coupled oscillator network can not achieve
synchronization.\hfill$\square$

\begin{lemma}\quad 
If there exists a $t_0\geq0$ such that
$\varphi_1(t_0)=\cdots=\varphi_m(t_0)=\varphi\quad(0<\varphi<1)$,
$D_1(t_0)=\cdots=D_m(t_0)=\{\eta\}\quad(0<\eta<\tau)$ and
$\varphi_{m+1}(t_0)=\cdots=\varphi_N(t_0)=1$,
$D_{m+1}(t_0)=\cdots=D_N(t_0)=\emptyset$, then the pulse-coupled
oscillator network can not achieve synchronization.
\end{lemma}
{\bf Proof:}  Similarly as in the proof of Lemma 6, we only need
to discuss two oscillator $i$ and $j$, with $1\leq i \leq m$ and
$m+1\leq j \leq N$.

Let $n=N-m$.

By $D_1(t_0)=\cdots=D_{i-1}(t_0)=D_{i+1}(t_0)=
\cdots=D_m(t_0)=\{\eta\}\,\,(0<\eta<\tau)$ and Lemma 1, oscillator
$1,\ldots,i-1,i+1,\ldots,m$ do not fire in time
$[\max(0,t_0-(\tau-\eta)-\tau),\ \max(0,t_0-\tau)]$. By
$\varphi_{m+1}(t_0)=\cdots=\varphi_N(t_0)=1$ and Lemma 1, oscillator
$m+1,\ldots,N$ do not fire in $[\max(0,t_0-(\tau-\eta)-\tau),\
\max(0,t_0-\tau)]$, too. Therefore, no spikes can reach the $i$th
oscillator in time $[t_0-(\tau-\eta),t_0]$. Then we have
$\varphi_i(t_0)=\tau-\eta$, i.e., $\eta=\tau-\varphi$.

Now we have the following table of $\varphi_i(t)$, $D_i(t)$,
$\varphi_j(t)$ and $D_j(t)$:
{\footnotesize
\begin{center}
\begin{tabular}{|l|l|l|l|l|}\hline
time $t$ & $\varphi_i(t)$ & $D_i(t)$ & $\varphi_j(t)$ &
$D_j(t)$\\
\hline \hline
$t_0$ & $\varphi$ & $\{\tau-\varphi\}$ & 1 & $\emptyset$\\
\hline
$t_0^{+}$ & $\varphi$ & $\{\tau-\varphi\}$ & 0 & $\{\tau\}$\\
\hline
$t_0+\tau-\varphi$ & $F_{m-1}(\tau)$ & $\{0\}$ & $F_m(\tau-\varphi)$ & $\{\varphi\}$\\
\hline
$(t_0+\tau-\varphi)^{+}$ & $F_{m-1}(\tau)$ & $\emptyset$ & $F_m(\tau-\varphi)$ & $\{\varphi\}$\\
\hline
$t_0+\tau$ & $F_n(F_{m-1}(\tau)+\varphi)$ & $\emptyset$ & $F_{n-1}(F_m(\tau-\varphi)+\varphi)$ & $\{0\}$\\
\hline
$(t_0+\tau)^{+}$ & $F_n(F_{m-1}(\tau)+\varphi)$ & $\emptyset$ & $F_{n-1}(F_m(\tau-\varphi)+\varphi)$ & $\emptyset$\\
\hline
\end{tabular}
\end{center}}

\bigskip
Denote $\phi_1=F_n(F_{m-1}(\tau)+\varphi)$ and
$\phi_2=F_{n-1}(F_m(\tau-\varphi)+\varphi)$.

Proposition 2 shows that
$\phi_1=F_n(F_{m-1}(\tau)+\varphi)>F_{m+n-1}(\tau)>
F_{n-1}(F_m(\tau-\varphi)+\varphi)=\phi_2$. Then
$\varphi_i(t_0+\tau+1-\phi_1)=1$,
$D_i(t_0+\tau+1-\phi_1)=\emptyset$,
$\varphi_j(t_0+\tau+1-\phi_1)=\phi_2+1-\phi_1$ and
$D_j(t_0+\tau+1-\phi_1)=\emptyset$. Thus, by Lemma 6, the
pulse-coupled oscillator network can not achieve synchronization.
\hfill$\square$

\section{Main theorem}\quad
By using the lemmas above, we can obtain our main theorem.

\begin{theorem}\quad 
Under the assumptions {\upshape(A1)} and {\upshape(A2)}, from any
initial phases {\upshape(}other than
$\varphi_1(0)=\cdots=\varphi_N(0)${\upshape)}, the pulse-coupled
oscillator network with delayed excitatory coupling can not achieve
complete synchronization.
\end{theorem}
{\bf Proof:} If $N=2$, then by Lemma 6 the theorem has been
proved. Therefore, we only need to consider the case of $N\geq3$.
When $N\geq3$, for any initial phases (other than
$\varphi_1(0)=\cdots=\varphi_N(0)$), one of the following two
cases must hold.

{\bf Case 1}: There exists a $t_1\geq0$ such that the network is
divided into two cliques (A and B) at time $t_1$. That is, the
oscillators of clique A have been synchronized at $t_1$ and the
oscillators of clique B have also been synchronized at $t_1$, but
the whole network has not been synchronized at $t_1$.

Because the initial phases do not satisfy
$\varphi_1(0)=\cdots=\varphi_N(0)$, there must be a $t_0$ with
$0\leq t_0\leq t_1$ such that the network is divided into two
cliques (A and B) at time $t_0$, and the network is divided into
three or more cliques at $t'$ for all $0\leq t'<t_0$. If $t_0=0$,
then the network is divided into two cliques at time 0. By Lemma 6,
it can not achieve synchronization. If $t_0>0$, similarly as in the
proof of Lemma 5, there must be an oscillator $i$ satisfying
$D_i(t_0)=\{0\}$. Then by Lemma 2(c), we have $0<\varphi_i(t_0)<1$.
Since the order of the oscillators is inessential to
synchronization, we can re-arrange the order and assume that clique
A includes oscillator $1,\ldots,m$ and clique B includes oscillator
$m+1,\ldots,N$, with $1\leq m<N$. Without loss of generality, we
assume that oscillator $i$ is in clique A, i.e., $1\leq i\leq m$. By
Lemma 3, we have $0<\varphi_1(t_0)=\cdots=\varphi_m(t_0)<1$ and
$D_1(t_0^{+})=\cdots=D_m(t_0^{+})=\emptyset$. For clique B, one of
the following two cases must hold.

{\bf Subcase 1.1}: There does not exist $t''<t_0$ such that the
oscillators of clique B have been synchronized at $t''$.

By Lemma 4, $D_i(t_0)=\{0\}$ shows $D_j(t_0)\neq\{0\}$ for all
$m+1\leq j \leq N$. Then by Lemma 5, we get
\begin{eqnarray*}
&&\varphi_1(t_0)=\cdots=\varphi_m(t_0)<1\\
&&D_1(t_0^{+})=\cdots=D_m(t_0^{+})=\emptyset\\
&&\varphi_{m+1}(t_0)=\cdots=\varphi_N(t_0)=1\\
&&D_{m+1}(t_0)=\cdots=D_N(t_0)=\emptyset.
\end{eqnarray*}
So
\begin{eqnarray*}
&&\varphi_1(t_0^{+})=\cdots=\varphi_m(t_0^{+})<1\\
&&D_1(t_0^{+})=\cdots=D_m(t_0^{+})=\emptyset\\
&&\varphi_{m+1}(t_0^{+})=\cdots=\varphi_N(t_0^{+})=0\\
&&D_{m+1}(t_0^{+})=\cdots=D_N(t_0^{+})=\{\tau\}.
\end{eqnarray*}
Similarly as in the proof of Lemma 6, the network can not achieve
synchronization.

{\bf Subcase 1.2}: There exists a $t''<t_0$ such that the
oscillators of clique B have been synchronized at $t''$.

By Lemma 4, we have $D_{m+1}(t_0)=\cdots=D_N(t_0)$. There are two
cases again.

(i) $D_{m+1}(t_0)=\cdots=D_N(t_0)=\emptyset$.

Then we have
\begin{eqnarray*}
&&\varphi_1(t_0)=\cdots=\varphi_m(t_0)<1\\
&&D_1(t_0^{+})=\cdots=D_m(t_0^{+})=\emptyset\\
&&\varphi_{m+1}(t_0)=\cdots=\varphi_N(t_0)\leq1\\
&&D_{m+1}(t_0)=\cdots=D_N(t_0)=\emptyset
\end{eqnarray*}

Denote $\varphi_1(t_0)=\cdots=\varphi_m(t_0)=\phi_1$ and
$\varphi_{m+1}(t_0)=\cdots=\varphi_N(t_0)=\phi_2$. If
$\phi_1=\phi_2$, then
\begin{eqnarray*}
&&\varphi_1(t_0)=\cdots=\varphi_m(t_0)=\varphi_{m+1}(t_0)=\cdots=\varphi_N(t_0)=\phi_1<1\\
&&D_1(t_0^{+})=\cdots=D_m(t_0^{+})=D_{m+1}(t_0^{+})=\cdots=D_N(t_0^{+})=\emptyset.
\end{eqnarray*}
By Lemma 3, the network has been synchronized at $t_0$, which
contradicts that $t_0\leq t_1$. Thus, $\phi_1\neq\phi_2$. Without
loss of generality, we may assume $\phi_1<\phi_2$. So, we have
\begin{eqnarray*}
&&\varphi_1(t_0+1-\phi_2)=\cdots=\varphi_m(t_0+1-\phi_2)=\phi_1+1-\phi_2<1\\
&&D_1(t_0+1-\phi_2)=\cdots=D_m(t_0+1-\phi_2)=\emptyset\\
&&\varphi_{m+1}(t_0+1-\phi_2)=\cdots=\varphi_N(t_0+1-\phi_2)=1\\
&&D_{m+1}(t_0+1-\phi_2)=\cdots=D_N(t_0+1-\phi_2)=\emptyset
\end{eqnarray*}
By Lemma 6, the network can not achieve synchronization.

(ii) $D_{m+1}(t_0)=\cdots=D_N(t_0)=\{\eta\}$ ($D_i(t_0)=\{0\}$
shows $0<\eta<\tau$).

Then we have
\begin{eqnarray*}
&&\varphi_1(t_0)=\cdots=\varphi_m(t_0)<1\\
&&D_1(t_0^{+})=\cdots=D_m(t_0^{+})=\emptyset\\
&&\varphi_{m+1}(t_0)=\cdots=\varphi_N(t_0)<1\\
&&D_{m+1}(t_0)=\cdots=D_N(t_0)=\{\eta\}
\end{eqnarray*}

Denote $\varphi_1(t_0)=\cdots=\varphi_m(t_0)=\phi_1$ and
$\varphi_{m+1}(t_0)=\cdots=\varphi_N(t_0)=\phi_2$. It follows from
Lemma 1 that $\phi_2+\eta<1$.

If $\phi_1+\eta>1$, then $\phi_1>\phi_2$. It follows that
\begin{eqnarray*}
&&\varphi_1(t_0+1-\phi_1)=\cdots=\varphi_m(t_0+1-\phi_1)=1\\
&&D_1(t_0+1-\phi_1)=\cdots=D_m(t_0+1-\phi_1)=\emptyset\\
&&\varphi_{m+1}(t_0+1-\phi_1)=\cdots=\varphi_N(t_0+1-\phi_1)=\phi_2+1-\phi_1<1\\
&&D_{m+1}(t_0+1-\phi_1)=\cdots=D_N(t_0+1-\phi_1)=\{\eta-(1-\phi_1)\}
\end{eqnarray*}
where $0<\eta-(1-\phi_1)<\tau$. By Lemma 7, the network can not
achieve synchronization.

If $\phi_1+\eta\leq1$, then $\phi_1\leq\phi_2$. It follows that
\begin{eqnarray*}
&&\varphi_1(t_0+\eta)=\cdots=\varphi_m(t_0+\eta)=F_{N-m}(\phi_1+\eta)\\
&&D_1(t_0+\eta)=\cdots=D_m(t_0+\eta)=\emptyset\\
&&\varphi_{m+1}(t_0+\eta)=\cdots=\varphi_N(t_0+\eta)=F_{N-m-1}(\phi_2+\eta)\\
&&D_{m+1}(t_0+\eta)=\cdots=D_N(t_0+\eta)=\{0\}.
\end{eqnarray*}
Let $K$ be a set of oscillators: oscillator $k'\in K$ if and only if
$D_{k'}(t_0)=\{0\}$. And let $n$ be the number of the oscillators in
$K$. It is obvious that clique B $\cap K=\emptyset$ and clique A
$\supseteq K$. From Proposition 1 and $\phi_1,\phi_2<1$, it follows
that
$\phi_1=\varphi_i(t_0)=f^{-1}(f(\varphi_i(t_0^{-}))+(n-1)\varepsilon)$
and
$\phi_2=\varphi_{m+1}(t_0)=f^{-1}(f(\varphi_{m+1}(t_0^{-}))+n\varepsilon)$.
Since $D_i(t_0)=\{0\}$ and $D_{m+1}(t_0)=\{\eta\}$ ($0<\eta<\tau$),
by Lemma 2(d), we have $\varphi_i(t_0^{-})>\varphi_{m+1}(t_0^{-})$.
So, $f(\phi_2)-f(\phi_1)=[f(\varphi_{m+1}(t_0^{-}))+n\varepsilon]-
[f(\varphi_i(t_0^{-}))+(n-1)\varepsilon]<\varepsilon$. From the
concavity assumption of function $f$, it follows that
$f(\phi_2+\eta)-f(\phi_1+\eta)<\varepsilon$. Thus,
$F_{N-m}(\phi_1+\eta)\neq F_{N-m-1}(\phi_2+\eta)$. Then, similarly
as in the proof of Lemma 6, the network can not achieve
synchronization.

{\bf Case 2}: There does not exist $t_1\geq0$ such that the network
is divided into two cliques at time $t_1$.

Assume that the pulse-coupled network can achieve synchronization.

Then there must be a $t_0>0$ such that the network has been
synchronized at $t_0$, but the network has not been synchronized at
$t$ for all $0\leq t<t_0$. Similarly as in the proof of Lemma 5,
there must be an oscillator $i$ such that $D_i(t_0)=\{0\}$. Because
there does not exist $t_1\geq0$ such that the network is divided
into two cliques at time $t_1$, there must be two other oscillators
$j$ and $k$ such that
\begin{eqnarray*}
&&\mbox{oscillator $i$ and $j$ have not been synchronized at $t$ for all $t<t_0$};\\
&&\mbox{oscillator $j$ and $k$ have not been synchronized at $t$ for all $t<t_0$};\\
&&\mbox{oscillator $k$ and $i$ have not been synchronized at $t$ for
all $t<t_0$}.
\end{eqnarray*}
By Lemma 4, $D_i(t_0)=\{0\}$ shows that $D_j(t_0)\neq\{0\}$ and
$D_k(t_0)\neq\{0\}$. So, by Lemma 5, we have
\begin{eqnarray*}
&&\varphi_j(t_0^{-})\neq\varphi_k(t_0^{-})\\
&&\varphi_j(t_0)=\varphi_k(t_0)=1\\
&&D_j(t_0)=D_k(t_0)=\emptyset
\end{eqnarray*}
However, $\varphi_i(t_0)=\varphi_j(t_0)=\varphi_k(t_0)=1$
contradicts $D_i(t_0)=\{0\}$ by Lemma 2(c).

So, in this case, the network can not achieve synchronization too.

In summary, the theorem is proved.\hfill$\square$

Theorem 1 shows that the presence of transmission delays can lead to
desynchronization. However, one would say that the definition of
complete synchronization in the theorem is awfully strong: it
requires that all oscillators reach perfect synchronization in
finite time; and say that usually synchronization only requires that
the phase differences between oscillators converge to zero as time
goes to infinity, i.e., asymptotical synchronization. In Mirollo and
Strogatz's model this distinction is irrelevant, since if the phase
differences are sufficiently small, perfect synchronization will be
achieved after the next round of firings. But this is no longer the
case when there exist transmission delays. Large numbers of
numerical examples show that when there exist transmission delays,
the phase differences can be small but there could still be big
variations in the times to process the pulses in the pipeline.
Simulation 1 in the following gives an example. Basing on the
numerical result, we conjecture that asymptotical synchronization
also cannot be achieved in pulse-coupled oscillator networks with
delayed excitatory coupling.

If our conjecture is true, then another problem arises: since
neither complete synchronization nor asymptotical synchronization
can be achieved, what does the model do? We numerically found that
the model usually converges to two or more phased-locked clusters.
It may be used to explain the phenomenon of clustering
synchronization in some species of fireflies, which follow the phase
advance model \cite{Buck1988}. Simulation 2 in the following shows
this.

\section{Numerical simulations}\quad\,
In this section, we give two numerical simulations.

Consider the pulse-coupled oscillator network with delayed
excitatory coupling:
\begin{eqnarray*}
&&f(\varphi)=-I\,\exp\{\ln\frac{I-1}{I}\cdot\varphi\}+I\qquad(I=1.05)\\
&&N=100\\
&&\varepsilon=0.001\\
&&\tau=0.1
\end{eqnarray*}
Obviously, it satisfies all the conditions of Theorem 1. In
particular, we have $f(2\tau)+N\varepsilon=0.5789<1$. So, we
conclude that the network can not achieve complete synchronization.

\subsection*{\itshape\mdseries{Simulation 1}}\quad
In this simulation, we illustrate that even if all the initial
phases are very close to each other, there could still be big
variations in the times to process the pulses in the pipeline. An
example with the initial phases chosen from a uniform distribution
on $(0,0.01]$ is plotted in Fig.4.

\setlength{\unitlength}{1cm}
\begin{picture}(9,9.5)
\put(2,0.5){\epsfig{file=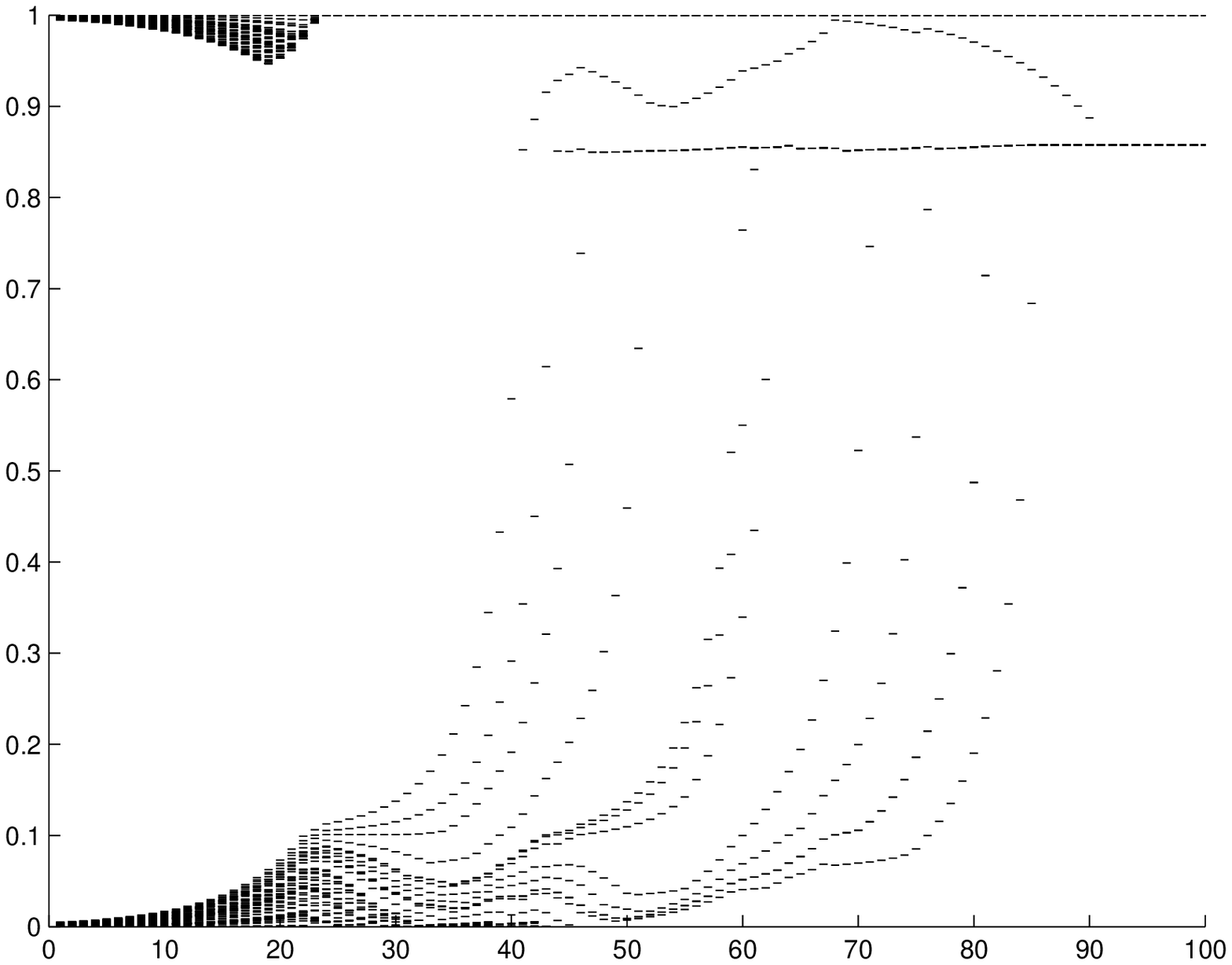,height=8cm}}
\put(1.5,4.5){\rotatebox{90}{\makebox(0,0)[c]{\footnotesize$\varphi_i(t_k)$}}}
\put(7.2,0.2){\makebox(0,0)[c]{\footnotesize$k$}}
\end{picture}

{\footnotesize Fig.4: Stroboscopic view on the phases
$\varphi_i(t_k)$ of $N=100$ oscillators plotted each time $t_k$ the
first oscillator fires its $k$th time, i.e., $\varphi_1(t_k)=1$. The
initial phases are chosen from a uniform distribution on
$(0,0.01]$.}

\subsection*{\itshape\mdseries{Simulation 2}}\quad
In this simulation, we investigate the clustering phenomena. The
oscillators are initialized with a uniform random distribution on
$(0,1]$. Fig.5 shows that the network can converge to two or more
phased-locked clusters.

\setlength{\unitlength}{1cm}
\begin{picture}(9,9)
\put(0,8){(a)}
\put(2,0.5){\epsfig{file=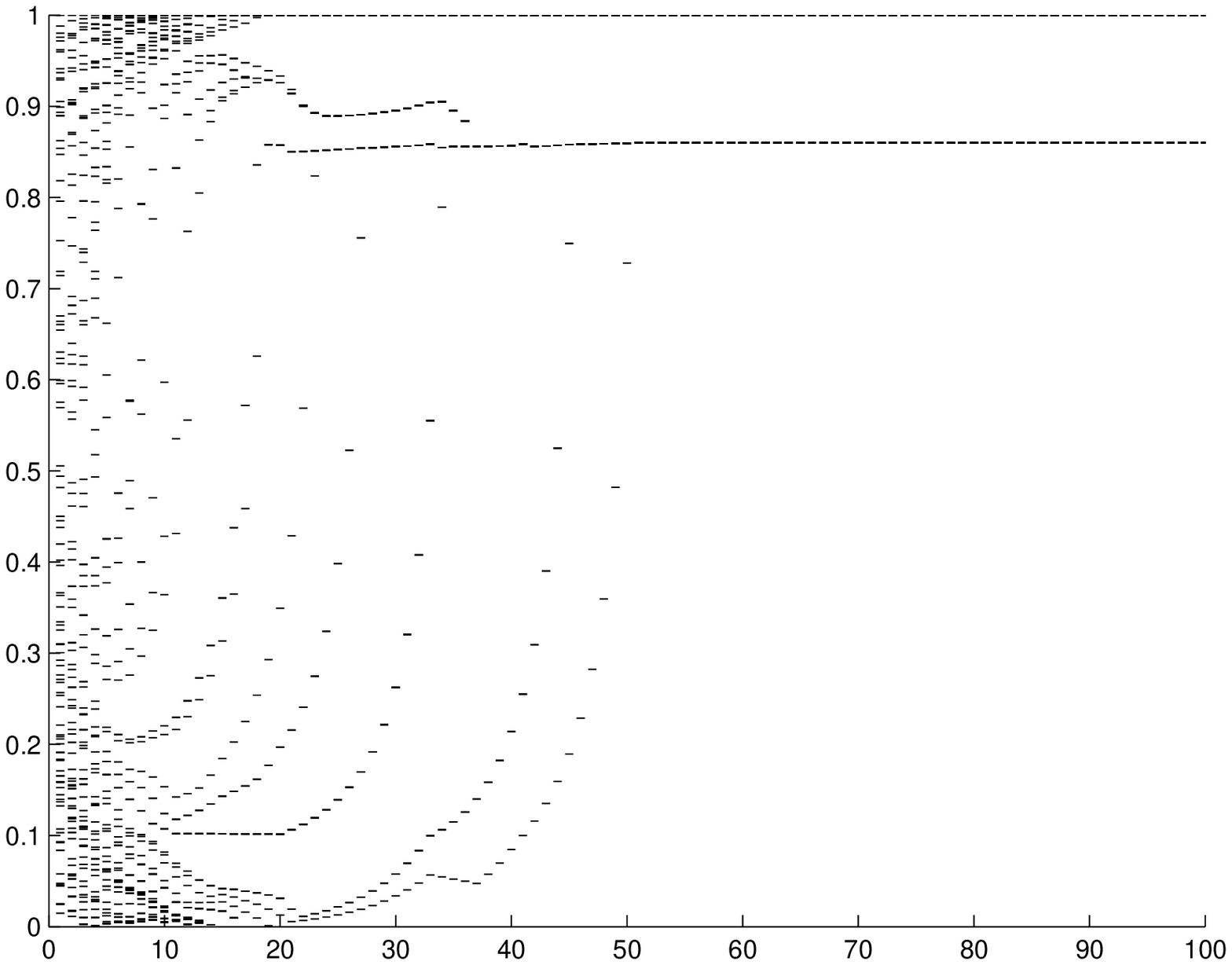,height=8cm}}
\put(1.5,4.5){\rotatebox{90}{\makebox(0,0)[c]{\footnotesize$\varphi_i(t_k)$}}}
\put(7.2,0.2){\makebox(0,0)[c]{\footnotesize$k$}}
\put(12.4,8){\makebox(0,0)[c]{\small Cluster 1}}
\put(11.5,8){\vector(-2,1){0.4}}
\put(12.4,6.9){\makebox(0,0)[c]{\small Cluster 2}}
\put(11.5,6.9){\vector(-2,1){0.4}}
\end{picture}

\setlength{\unitlength}{1cm}
\begin{picture}(9,8.8)
\put(0,8){(b)}
\put(2,0.5){\epsfig{file=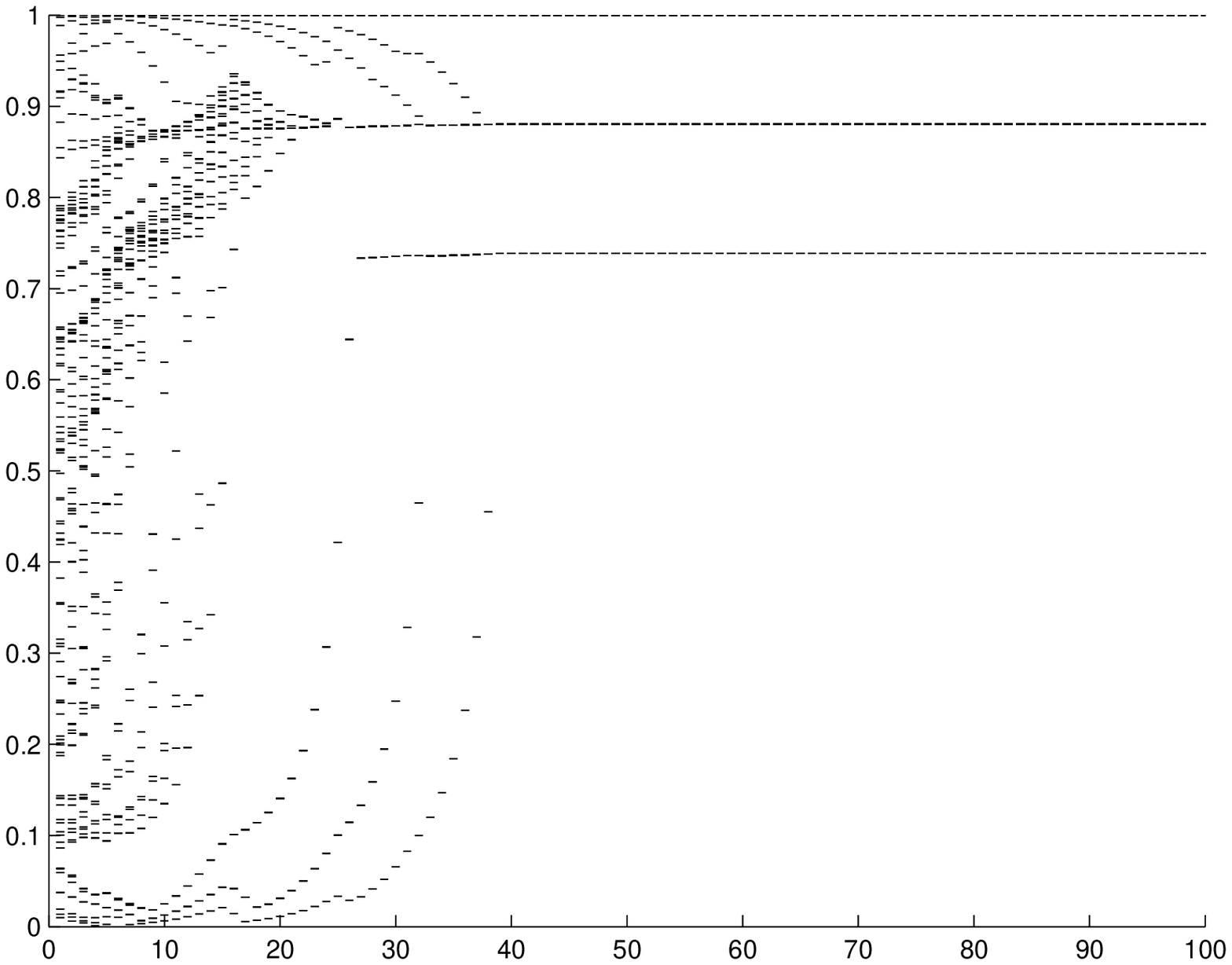,height=8cm}}
\put(1.5,4.5){\rotatebox{90}{\makebox(0,0)[c]{\footnotesize$\varphi_i(t_k)$}}}
\put(7.2,0.2){\makebox(0,0)[c]{\footnotesize$k$}}
\put(12.4,8){\makebox(0,0)[c]{\small Cluster 1}}
\put(11.5,8){\vector(-2,1){0.4}}
\put(12.4,7.1){\makebox(0,0)[c]{\small Cluster 2}}
\put(11.5,7.1){\vector(-2,1){0.4}}
\put(12.4,6){\makebox(0,0)[c]{\small Cluster 3}}
\put(11.5,6){\vector(-2,1){0.4}}
\end{picture}

{\footnotesize Fig.5: Stroboscopic view on the phases
$\varphi_i(t_k)$ of $N=100$ oscillators plotted each time $t_k$ the
first oscillator fires its $k$th time, i.e., $\varphi_1(t_k)=1$. The
initial phases are chosen from a uniform distribution on $(0,1]$.
Eventually, the network is divided into two or more phased-locked
clusters, for example, (a) two clusters (b) three clusters.}

\section{Conclusions}\quad
In this paper, pulse-coupled oscillator networks with delayed
excitatory coupling are studied. We propose an assumption (A2),
which is reasonable to real biological systems, especially to
fireflies. It is proved that under the assumptions (A1) and (A2),
from any initial phases (other than
$\varphi_1(0)=\cdots=\varphi_N(0)$), the network can not achieve
complete synchronization.  This result can explain why {\itshape
Photinus pyralis} rarely synchronizes flashing, which is known as an
example of pulse-coupled oscillator networks with delayed excitatory
coupling. Furthermore, according to Simulation 1, we conjecture that
asymptotical synchronization also cannot be achieved; and in
Simulation 2, we exhibit a phenomenon of clustering synchronization.


\begin{thebibliography}{99}
\bibitem{Buck1976} 
J. Buck, E. Buck, Synchronous fireflies, Sci. Am. 234 (1976)
74-85.

\bibitem{Buck1988} 
J. Buck, Synchronous rhythmic flashing of fireflies, Part II. Q.
Rev. Biol. 63 (3) (1988) 265-289.

\bibitem{Hanson1978} 
F.E. Hanson, Comparative studies of firefly pacemakers, Fed.
Proc., 37 (1978) 2158-2164.

\bibitem{Smith1935} 
H.M. Smith, Synchronous flashing of fireflies, Science, 82 (1935)
151.

\bibitem{Peskin1975} 
C.S. Peskin, Mathematical aspects of heart physiology, Courant
Institute of Mathematical Sciences, New York University, New York,
(1975) 268-278.

\bibitem{Mir1990} 
R.E. Mirollo and S.H. Strogatz, Synchronization of pulse-coupled
biological oscillators, SIAM J. Appl. Math., 50 (1990) 1645-1662.

\bibitem{Kuramoto1991} 
Y. Kuramoto, Collective synchronization of pulse-coupled
oscillators and excitable units, Phys. D., 50 (1991) 15-30.

\bibitem{Vanvreeswijk1993} 
C. Vanvreeswijk, L.F. Abbott, Self-sustained firing in populations
of integrate-and-fire neurons, SIAM J. Appl. Math., 53 (1993)
253-264.

\bibitem{Goel2002} 
P. Goel, B. Ermentrout, Synchrony, stability, and firing patterns
in pulse-coupled oscillators, Phys. D., 163 (2002) 191-216.

\bibitem{Chen1994} 
C.C. Chen, Threshold effects on synchronization of pulse-coupled
oscillators, Phys. Rev. E, 49 (1994) 2668-2672.

\bibitem{Corral1995} 
A. Corral, C.J. Perez, A. Diazguilera, et al, Self-organized
criticality and synchronization in a lattice model of
integrate-and-fire oscillators, Phys. Rev. Lett., 74 (1995)
118-121.

\bibitem{Mathar1996} 
R. Mathar, J. Mattfeldt, Pulse-coupled decentral Synchronization,
SIAM J. Appl. Math., 56 (1996) 1094-1106.

\bibitem{Nischwitz1995} 
A. Nischwitz, H. Glunder, ``Local lateral inhibition-a key to
spike synchronization'', Biol. Cybern., 73 (1995) 389-400.

\bibitem{Ernst1995} 
U. Ernst, K. Pawelzik, T. Geisel, Synchronization induced by
temporal delays in pulse-coupled oscillators, Phys. Rev. Lett., 74
(1995) 1570-1573.

\bibitem{Ernst1998} 
U. Ernst, K. Pawelzik, T. Geisel, Delay-induced multistable
synchronization of biological oscillators, Phys. Rev. E., 57
(1998) 2150-2162.

\bibitem{Knoblauch2002} 
A. Knoblauch, G. Palm, ``Scene segmentation by spike
synchronization in reciprocally connected visual areas. II. Global
assemblies and synchronization on larger space and time scales'',
Biol. Cybern., 87 (2002) 168-184.

\bibitem{Coombes1997} 
S. Coombes, G.J. Lord, Desynchronization of pulse-coupled
integrate-and-fire neurons, Phys. Rev. E., 55 (1997) 2104-2107.

\bibitem{Kim2004} 
D.E. Kim, A spiking neuron model for synchronous flashing of
fireflies, Biosystems, 76 (2004) 7-20.
\end{thebibliography}
\end{document}